\documentclass[runningheads]{llncs}
\pdfoutput=1
\usepackage{graphicx}
\usepackage{amsmath}
\usepackage{amssymb}
\usepackage{todonotes}
\usepackage{bbm} 
\usepackage{longtable}
\usepackage{aliascnt} 
\usepackage{hyperref}
\usepackage{multirow} 
\usepackage{enumitem} 
\setenumerate[0]{label=(\alph*)}  
\usepackage{algorithm}
\usepackage{algpseudocode}
\floatname{algorithm}{Algorithm}

\DeclareMathOperator{\STAB}{STAB}
\DeclareMathOperator{\COL}{COL}
\DeclareMathOperator{\CUT}{CUT}
\DeclareMathOperator{\conv}{conv}
\DeclareMathOperator{\trace}{trace}
\DeclareMathOperator{\diag}{diag}

\newcommand{\COLE}{\COL^{\varepsilon}}
\newcommand{\eshmc}[1]{z_{mc}^{#1}}
\newcommand{\eshss}[1]{z_{ss}^{#1}}
\newcommand{\eshc}[1]{z_{c}^{#1}}
\newcommand{\eshcE}[1]{z_{c\varepsilon}^{#1}}

\newcommand{\XI}{X_{I}}
\newcommand{\GI}{G_{I}}
\newcommand{\bI}{b_{I}}
\newcommand{\kI}{k_{I}}
\newcommand{\yI}{y_{I}}
\newcommand{\yJ}{y}
\newcommand{\lambdaI}{\lambda_{I}}
\newcommand{\lambdaJ}{\lambda}

\newcommand{\opvt}{ {\transposedVec} }

\newcommand{\sumbI}{b}
\newcommand{\tI}{t_{I}}
\newcommand{\setX}{{\mathcal X}}
\newcommand{\setXmc}{{\mathcal X}^{E}}
\newcommand{\setXss}{{\mathcal X}^{S}}
\newcommand{\setXcol}{{\mathcal X}^{C}}

\newcommand{\R}{{\mathbb R}}
\newcommand{\Sym}[1]{{\mathcal S}_{#1}}
\newcommand{\Lagrangian}{{\mathcal L}}
\newcommand{\simplex}{{\Delta}}
\newcommand{\allones}[1]{\mathbbm{1}_{#1}}
\newcommand{\transposedVec}{T}
\newcommand{\transposedOp}{\top}
\newcommand{\norm}[1]{\left\lVert #1 \right\rVert}

\newcommand{\opA}{{\mathcal A}_I}
\newcommand{\opAt}{{\mathcal A}_I^{\transposedOp}}
\newcommand{\opAmc}{{\mathcal A}_I}
\newcommand{\opAmct}{{\mathcal A}_I^{\transposedOp}} 
\newcommand{\opM}{{\mathcal M}_I}
\newcommand{\opMt}{{\mathcal M}_I^{\transposedOp}}
\newcommand{\opMmc}{{\mathcal M}_I}
\newcommand{\opMmct}{{\mathcal M}_I^{\transposedOp}} 
\newcommand{\opP}{{\mathcal P}_I}
\newcommand{\opPt}{{\mathcal P}_I^{\transposedOp}}
\newcommand{\opD}{{\mathcal D}_I}
\newcommand{\opDt}{{\mathcal D}_I^{\transposedOp}}

\newcommand{\hBdlExplanation}[1]{h_{#1}}
\newcommand{\yBdlExplanation}[1]{y_{#1}}
\newcommand{\gBdlExplanation}[1]{g_{#1}}
\newcommand{\eBdlExplanation}[1]{e_{#1}}
\newcommand{\XBdlExplanation}[1]{X_{#1}}

\newcommand{\yCenter}{\overline{y}}
\newcommand{\yTrial}{\widetilde{y}}

\newcommand{\hCenter}{\overline{h}}
\newcommand{\wTrial}{\widetilde{w}}
\newcommand{\vTrial}{\widetilde{v}}
\newcommand{\vI}{v_{I}}
\newcommand{\bdlSize}{r}
\newcommand{\bdlPar}{\mu}
\newcommand{\nrESC}{q}
\newcommand{\Bundle}{\mathcal{B}}

\DeclareMathOperator{\st}{s.t.}
\newcommand{\yICenter}{\overline{y}_{I}}
\newcommand{\opG}{{\mathcal G}}

\newcommand{\opGI}{{\mathcal G}_I}
\newcommand{\opGItrans}{{\mathcal G}_I^{\transposedOp}}

\newcommand{\betaI}{\beta_{I}}

\usepackage{colortbl}

\begin{document}
\title{A Computational Study of\\ Exact Subgraph Based SDP Bounds for Max-Cut, 
Stable Set and 
Coloring%
\thanks{This project has received funding from the European 
Union's Horizon~2020 research 
and innovation programme under the Marie Sk\l{}odowska-Curie grant 
agreement No~764759 and the 
Austrian Science Fund (FWF): I 3199-N31 and P 28008-N35.
We thank three anonymous referees for their constructive comments 
which 
substantially helped to improve the presentation of our material.
An extended abstract containing the theoretical foundations of this work 
appeared 
as~\cite{GaarRendl}.
This article now contains more details, additional information and many
new computational results.
}}

%
\titlerunning{Exact Subgraph Based SDP Bounds}
%
\author{Elisabeth Gaar\inst{1}\orcidID{0000-0002-1643-6066} \and
Franz Rendl\inst{1}\orcidID{0000-0003-1578-9414}}

\authorrunning{E. Gaar and F. Rendl}
%
\institute{Alpen-Adria-Universit\"{a}t Klagenfurt, Institut f\"{u}r 
Mathematik,\\
    Universit\"{a}tsstr. 65-67, 9020 Klagenfurt, Austria \\
    \email{\{elisabeth.gaar,franz.rendl\}@aau.at}}
\maketitle              
\begin{abstract}
The ``exact subgraph'' approach was recently introduced as a hierarchical
scheme to get increasingly tight semidefinite programming 
relaxations of several NP-hard 
graph optimization problems. 
Solving these relaxations is a computational challenge because 
of the potentially large number of violated subgraph constraints.
We introduce a computational framework for these relaxations designed
to cope with these difficulties. We suggest a partial 
Lagrangian dual, and exploit the fact that its
evaluation  decomposes into several independent subproblems.  
This opens the way to use the bundle method from non-smooth 
optimization to minimize the dual function.
Finally
computational experiments on the Max-Cut, stable set and coloring
problem show the excellent quality of the bounds obtained with this approach.

\keywords{Semidefinite programming  \and Relaxation hierarchy \and Max-Cut \and Stable set 
\and Coloring.}
\end{abstract}
%
%
%


\section{Introduction}
The study of NP-hard problems has led to the introduction of 
various hierarchies of relaxations, which typically involve several 
levels. Moving from one level to the next the relaxations get 
increasingly tighter and ultimately the exact optimum may be reached,
but the computational effort grows accordingly.

Among the most prominent hierarchies are the polyhedral ones from Boros, Crama and Hammer 
\cite{BorosCramaHammer} as well as the ones from Sherali and Adams \cite{SheraliAdamsHierarchy}, 
Lov\'{a}sz and Schrijver \cite{LovaszSchrijverHierarchy} and Lasserre \cite{LasserreHierarchy} 
which are based on semidefinite programming (SDP). Even though on the starting level they have a 
simple SDP relaxation, 
already the first nontrivial level in the hierarchy requires the 
solution of SDPs in matrices of order $\binom{n}{2}$ and on level $k$ the matrix order is 
$n^{O(k)}$. 
Hence they are considered mainly as theoretical tools and from a practical point of view  these 
hierarchies 
are of  limited use.

Not all hierarchies are of this type. In \cite{BorosCramaHammer} a polyhedral
hierarchy for the Max-Cut problem is introduced which maintains 
$\binom{n}{2}$ variables at all levels, with a growing number of 
constraints. 
More recently, 
Adams, Anjos, Rendl 
and Wiegele~\cite{AARW} introduced a hierarchy of SDP relaxations
which act in the space of symmetric $n \times n$ matrices and at level 
$k$ of the hierarchy all submatrices of order $k$ have to be ``exact'' in a 
well-defined sense, i.e. they have to fulfill an \emph{exact subgraph constraint} (ESC).

It is the main purpose of this paper to describe an efficient way
to optimize over level $k$ of this hierarchy for small values of 
$k$, e.g. $k\leqslant 7$, and demonstrate the efficiency
of our approach for the Max-Cut, stable set
and coloring problem. 
These investigations 
were started in \cite{elli-diss,GaarRendl} and here we offer the full picture.

Maintaining $\binom{n}{k}$ possible ESCs in 
an SDP in matrices of order $n$ is computationally 
infeasible even for $k=2$ or $k=3$, because each ESC
creates roughly $\binom{k}{2}$ additional equality constraints and at most $2^k$ additional 
 variables.

We suggest
the following ideas to overcome this difficulty. First we proceed
iteratively, and in each iteration we include only
(a few hundred of) the most  violated ESCs. 
More importantly, we propose to solve the dual of the resulting
SDP. The structure of this SDP with ESCs
admits a reformulation of the dual in the form of a non-smooth
convex minimization problem with attractive features.
First, any dual solution
yields a valid bound for our relaxations, so it is not necessary 
to carry out the minimization to optimality.
Secondly, the dual function 
evaluation decomposes into two independent problems. The first one
is simply a sum of max-terms (one for each ESC), 
and the second one consists in solving a ``basic'' SDP, independent
of the ESCs. The optimizer for this second problem
also yields a subgradient of the objective function. With this 
information at hand we suggest to use the bundle method from non-smooth 
convex optimization. It 
provides an effective machinery to get close to a minimizer
in few iterations. 

As a result we are able to get near optimal solutions where 
all ESCs for  small values of $k$ ($k \leqslant 7$)
are satisfied up to a small error tolerance. Our computational
results demonstrate the practical potential of this approach. 

The paper is organized as follows. In Section~\ref{sec:DefProblemsBasicRel} we briefly describe the 
Max-Cut, the stable set and the coloring problem along with their
semidefinite relaxations, which are well-studied in the literature.
Section~\ref{sec: ESC} recalls the exact subgraph hierarchy, described in \cite{AARW}.
We introduce a unified setting for all these problems and take a look at
their structural properties. 
In Section~\ref{sec: lagrangianDual} we reformulate the SDP and consider a partial Lagrangian dual. 
It results in many subproblems, separating the basic SDP part from the ESC 
part. 
The bundle method from non-smooth optimization is described in Section~\ref{sec:bundle} as
an attractive algorithmic framework to deal with the subproblems in the
partial Lagrangian dual. 
In Section~\ref{sec:overallAlgo} we describe our algorithm in order to obtain 
exact subgraph based SDP bounds. 
We argue in Section~\ref{sec:comp} that standard SDP solvers are 
only of limited use when dealing with our ESC hierarchy and present 
extensive computational results.
Finally we close with conclusions and future work in  
Section~\ref{sec:conclusions}.

We finish this introductory section with some notation. 
We denote the vector of all-ones of size $n$ with $\allones{n}$ and $\simplex_n = \{x \in 
\R^{n}_{+}: \sum_{i=1}^{n}x_{i} = 1\}$. If the dimension is clear from the context we may omit 
the index and write $\allones{}$ and $\simplex$. Furthermore let $N = \{1, 2, \dots, n\}$. 
A graph $G$ on $n$ vertices has vertex set $N$ and edge set $E$.
The complement graph $\overline{G}$ of a graph $G$ has the same vertex 
set $N$ 
and contains an edge $\{i,j\} \subseteq N$ if and only if $\{i,j\} 
\not 
\in E$. 
$\Sym{n}$ 
is the set of $n$-dimensional symmetric matrices.
A spectrahedron is a set that is obtained as the intersection of the cone of 
positive semidefinite matrices with some linear affine subspace.


\section{Combinatorial Problems and Semidefinite Relaxations}
\label{sec:DefProblemsBasicRel} 
\subsection{The Max-Cut Problem}
In the Max-Cut problem a symmetric matrix $L \in
\Sym{n}$ is given and $c \in \{-1,1 \}^{n}$ which 
maximizes $c^{\transposedVec}Lc$ should be determined. 

If the matrix $L$ corresponds to the 
Laplacian matrix of a (edge-weighted undirected) graph $G$, 
this is equivalent to finding a partition of the vertices of $G$ 
into two subsets such
that  the total weight of the edges joining these two subsets
is maximized. Such an edge set is also called  a \emph{cut} in  $G$.

Partitions of $N$ into two subsets can be expressed as $c \in \{-1,1 \}^{n}$ 
where the 
two subsets of $N$ correspond to the entries of $c$ with the same sign. 
Given $c \in \{-1,1\}^{n}$ we call $C=cc^{\opvt}$ a \emph{cut
    matrix}. 
The convex hull of all cut matrices (of order $n$) is denoted by
$$\CUT_{n} = 
\conv \left\{ cc^{\opvt}:~ c \in \{-1,1\}^{n} \right\}
$$ 
or simply $\CUT$ if the dimension is clear from the context. 
Since $c^{\opvt}Lc = \langle L, cc^{\opvt} \rangle$  
the Max-Cut problem can also be written as the following (intractable) linear
program 
$$
z_{mc} = \max \{ \langle L, X\rangle:~ X \in \CUT \}.
$$
$\CUT$ is contained in the spectrahedron
$$
    \setXmc = \left\{ X \in \Sym{n} : \diag(X) = 
    \allones{n},\, X \succcurlyeq 0 \right\},
$$
hence
\begin{equation}
\label{relaxation mc}
r_{mc} = \max \left\{\langle L,X\rangle :~ X \in \setXmc \right\}
\end{equation}
is a basic semidefinite relaxation for Max-Cut. 
This model is well-known, attributed to Schrijver and  was introduced in a dual 
form by Delorme and Poljak~\cite{DelormePoljak}. 
It can be solved in polynomial time to a fixed prescribed precision and 
solving this relaxation for $n=1000$ takes only a few seconds.

It is well-known that the Max-Cut problem is NP-hard. 
On the positive side, Goemans and Williamson \cite{GoemansWilliamson}
show that one can find  a cut in a graph with
nonnegative edge weights of value at least 0.878$z_{mc}$ in polynomial time. 


\subsection{The Stable Set Problem}
In the stable set problem the input
is an unweighted graph $G$. 
We call a subset of the vertices \emph{stable}, if no two vertices are adjacent.
Moreover we call a vector $s \in \{0,1 \}^n$ a \emph{stable set vector} if it 
is the incidence vector of a stable set. 
The convex hull of all stable set
vectors of $G$  is denoted with $\STAB(G)$.  
In the stable set problem we want to determine the \emph{stability number $\alpha(G)$}, which 
denotes the cardinality of a largest stable set in $G$, hence
\begin{align*}
\alpha(G) = \max\left\{ \allones{}^{\opvt}s:~ s \in \STAB(G) \right\}.  
\end{align*}
Furthermore we denote with
\begin{equation*}
\STAB^{2}(G) =
\conv \left\{ ss^{\opvt}:~ s \in \STAB(G) \right\}
\end{equation*}
the convex hull of all \emph{stable set matrices} $ss^{\opvt}$. Then with the arguments of 
Gaar~\cite{elli-diss} it is easy to check that
\begin{align*}
\alpha(G) = \max\left\{ \trace(X):~ X \in \STAB^2(G) \right\}.
\end{align*}
Furthermore
$\STAB^{2}(G)$ is contained  
in the following spectrahedron 
\begin{equation*}
    \setXss = \left\{ X \in \Sym{n} :~ 
    X_{ij}=0 \quad 
    \forall \{i,j\} \in E,~ 
    x = \diag(X),~ 
    \left(
    \begin{array}{cc}
        1 & x^{\opvt} \\
        x & X
    \end{array}
    \right)
    \succcurlyeq 0 \right\},
\end{equation*}
which is known as the \emph{theta body} in the literature. Therefore 
\begin{equation}
\label{relaxation ss}
\vartheta(G)= \max \left\{ \trace(X):~ X \in \setXss \right\}
\end{equation}
is a relaxation of the stable set problem. The Lov{\'a}sz theta function $\vartheta(G)$ was 
introduced in a seminal paper by Lov{\'a}sz \cite{LovaszStart}.
We refer to Gr{\"o}tschel, Lov{\'a}sz and Schrijver 
\cite{OurUsedFormOfLovasTheta} for a comprehensive analysis of $\vartheta(G)$. 

Determining $\alpha(G)$ is again NP-hard. Contrary to Max-Cut, which
has a polynomial time .878-approximation, for every $\varepsilon>0$ there can be no polynomial 
time algorithm that 
approximates $\alpha(G)$ within a factor better than $O(n^{1-\varepsilon})$ unless
$P=NP$, see H{\aa}stad \cite{stableSetNotApproximable}.

\subsection{The Vertex Coloring Problem}
The coloring problem for a given graph $G$ consists in determining 
the \emph{chromatic number} $\chi(G)$, which is the smallest $t$ such that
$N$ can be partitioned into $t$ stable sets. 
Let $S=(s_{1}, \ldots, s_{k})$ be a matrix
where each column $s_i$ is a stable set vector 
and the corresponding stable sets 
partition
$N$ into $k$ sets.  
Let us call such  matrices $S$ {\em stable-set partition matrices} (SSPM)
and denote by $|S|$ the number of columns of $S$ or equivalently 
the number of stable set vectors of $S$.
The $n \times n$ matrix $X=SS^{T}$ is called {\em coloring matrix}. 
The convex hull of the
set of all coloring matrices of $G$
is denoted by
$$\COL(G) = \conv \left\{X:~ X \text{ is a coloring matrix of }G 
\right\}.$$
We also need the \emph{extended coloring polytope}
\begin{equation*}
\COLE(G) = \conv \left\{
\left( 
\begin{array}{cc} k & \allones{}^{\opvt}\\\allones{} & X
\end{array}\right)
= 
\sum_{i=1}^{k} 
\binom{1}{s_{i}}
\binom{1}{s_{i}}^{\opvt} 
:
\begin{array}{c}
S = (s_{1}, \ldots, s_{k}) \text{ is a} \\
\text{SSPM of } G,~ X = SS^{\opvt} 
\end{array}
\right\}.
\end{equation*}
The difficult set $\COLE$ can be relaxed to the easier spectrahedron
\begin{equation*} 
\setXcol = \left\{
\left( 
\begin{array}{cc} t & \allones{}^{\opvt}\\\allones{} & X 
\end{array}\right) \succcurlyeq 0:~ \diag(X)=\allones{n}, ~X_{ij}=0 ~\forall \{i,j\} \in E 
\right\}
\end{equation*}
and we can consider the semidefinite program  
\begin{equation}\label{relaxation col}
t^{*}(G) = \min \left\{t:~ 
\left( 
\begin{array}{cc} t & \allones{}^{\opvt}\\\allones{} & X 
\end{array}\right) \in \setXcol \right\}.
\end{equation}
Obviously 
$t^{*}(G) \leqslant \chi(G)$ holds
because the SSPM $S$ consisting of $\chi(G)$ stable sets 
yields a feasible coloring matrix $X=SS^{\opvt}$ with objective function value $\chi(G)$.
It is in fact a consequence of conic duality that
$t^{*}(G)= \vartheta(\overline{G})$ holds.


It is NP-hard to 
find  $\chi(G)$, to find a 4-coloring of a
3-colorable graph \cite{GuruswamiKhanna} and to color
a $k$-colorable graph with $O(k^{\frac{\log k}{25}})$ colors 
for sufficiently large $k$, \cite{Khot}.


\section{Exact Subgraph Hierarchy}\label{sec: ESC}
\subsection{Definition of the Hierarchy}
In this section we discuss how to systematically tighten the 
relaxations~\eqref{relaxation mc},~\eqref{relaxation ss} 
and~\eqref{relaxation col} with  
``exactness conditions'' imposed on small subgraphs. We obtained the 
relaxations 
by relaxing the 
feasible regions $\CUT$, 
$\STAB^{2}$ and $\COL$ of the integer problem to simple spectrahedral sets. Now 
we will use small 
subgraphs to get closer to the feasible regions of the original problems again.

For $I \subseteq N$ let $\kI=|I|$ be the cardinality of $I$.
Furthermore let $\GI$ be the induced subgraph of $G$ 
on the set of vertices $I$.
If $X$ is the $n \times n$ matrix from the relaxations~\eqref{relaxation 
mc},~\eqref{relaxation ss} or~\eqref{relaxation col}, 
then we denote with $\XI$ the principal $\kI \times \kI$ submatrix of $X$ 
corresponding to the rows 
and
columns in $I$. 
Note that $\XI$ is the submatrix of $X$ corresponding to $\GI$.

We first look at the exact subgraph relaxations for Max-Cut.
Adams, Anjos, Rendl and Wiegele~\cite{AARW} introduced additional  constraints 
for the Max-Cut relaxation~\eqref{relaxation mc} in the following way. 
The \emph{exact subgraph constraint} (ESC) for $I \subseteq N$
requires  
that the matrix $\XI$ corresponding to the subgraph $\GI$ lies in 
the convex hull of the cut matrices of $\GI$, 
that is
\begin{equation*}
\XI \in  \CUT_{|I|}. 
\end{equation*}
The ESC for $I$ can equivalently be phrased as  
\begin{equation*}
\XI = \sum_{i=1}^{\tI} \lambda_{i}C^{I}_{i}
\end{equation*}
for some $\lambda \in \simplex_{\tI}$ where $C^{I}_{i}$ is the 
$i$-th cut 
matrix of the subgraph $\GI$ and 
$\tI$ is the total number of cut matrices.
If $X$ is a solution of~\eqref{relaxation mc} that fulfills the ESC 
for some $I$ we say that $X$ is \emph{exact} on $I$ and 
$X$ is \emph{exact} on 
$\GI$.

Now we want the ESCs to be fulfilled not only 
for one but for a certain selection of  subgraphs.
We denote with $J$ the set of subsets $I$, on which we require $X$ to be exact,
and  get the following 
SDP relaxation with ESCs for Max-Cut.
\begin{equation}
\label{relaxation mc with esc}
\max \{\langle L,X\rangle:~
X \in \setXmc,~ \XI \in \CUT_{|I|} ~ \forall I \in J \}
\end{equation}

Before we give theoretical justification that~\eqref{relaxation mc with esc} is 
worth to be investigated, 
we present the ESCs for the other problems.
We start with the stable set problem on a graph $G$ and its 
relaxation~\eqref{relaxation ss}.
In this case the ESC for $I \subseteq N$, and hence for the subgraph $\GI$, 
requires 
that 
$
\XI \in \STAB^{2}(G_{I})
$
holds and 
the SDP with ESCs for the stable set problem is
\begin{equation}
\label{relaxation ss with esc}
\max \{ \trace(X):~ X \in \setXss,~
\XI \in \STAB^{2}(\GI) ~ \forall I \in J \}.
\end{equation}

Turning to  the coloring problem, we analogously impose additional ESCs
of the form
$
\XI \in \COL(G_I)
$
to obtain the SDP with ESCs
\begin{equation}
\label{relaxation col with esc}
\min \left\{t:~ 
\left( 
\begin{array}{cc} t & \allones{}^{\opvt}\\ \allones{} & X 
\end{array}\right) \in \setXcol,~
\XI \in \COL(\GI) ~ \forall I \in J
\right\}.
\end{equation}

We now want to investigate the properties 
of~\eqref{relaxation mc with esc},~\eqref{relaxation ss with esc} 
and~\eqref{relaxation col with esc}.
Towards that end we define the $k$-th level of the \emph{exact subgraph 
hierarchy} according to~\cite{AARW} 
by using $J = \{I \subseteq N:~ |I| = k\}$ in the 
SDPs~\eqref{relaxation mc with esc},~\eqref{relaxation ss with esc} 
and~\eqref{relaxation col with esc}, respectively. 
We denote the corresponding objective function 
values with $\eshmc{k}$, $\eshss{k}$ and $\eshc{k}$. So in other words the 
$k$-th level of the 
exact subgraph hierarchy is 
obtained by forcing all subgraphs on $k$ vertices to be exact in the basic SDP 
relaxation.

Note that 
\begin{alignat*}{10}
z_{mc} &= \eshmc{n}  &&\leqslant 
\dots 
&&\leqslant \eshmc{k} &&\leqslant \eshmc{k-1} &&\leqslant 
\dots
&&\leqslant \eshmc{2}
&&\leqslant \eshmc{1} &&= r_{mc}\\
\alpha(G) &= \eshss{n} &&\leqslant 
\dots 
&&\leqslant \eshss{k} &&\leqslant \eshss{k-1} &&\leqslant 
\dots
&&\leqslant \eshss{2} 
&&\leqslant \eshss{1} &&= \vartheta(G)
\end{alignat*}
holds for all $k \in \{2, \dots, n\}$, see~\cite{AARW,elli-diss}.
Hence~\eqref{relaxation mc with esc} and~\eqref{relaxation ss with esc} 
are relaxations of Max-Cut and the stable set problem.

Furthermore it can be verified that
\begin{alignat*}{10}
t^{*}(G) &= \eshc{1} &&\leqslant \eshc{2} &&\leqslant 
\dots 
&&\leqslant \eshc{k-1} &&\leqslant \eshc{k} &&\leqslant 
\dots
&&\leqslant \eshc{n} &&\leqslant \chi(G)
\end{alignat*}
holds for all $k \in \{2, \dots, n\}$, so for the coloring problem we do not
necessarily  
reach $\chi(G)$ at the $n$-th level. 
However, the following holds.
Let $\eshcE{k}$ be the optimal objective function value if we add the 
inequalities 
$t \geqslant \sum_{i=1}^{\tI}[\lambdaI]_{i}|S^{I}_{i}|$
where $\lambdaI \in \simplex_{\tI}$ is a 
variable for the 
convex combination for each subgraph $\GI$ to the SDP for $\eshc{k}$. Then 
$\eshcE{n} = \chi(G)$ 
holds. 
Hence $\eshc{k}$ is a relaxation of $\eshcE{k}$, which is in turn a 
relaxation of the coloring problem. As a result it is clear that it makes sense 
to investigate~\eqref{relaxation mc with esc},~\eqref{relaxation ss with esc} 
and~\eqref{relaxation col with esc}.

Note that in the case of the stable set and the coloring problem the polytopes $\STAB^{2}(\GI)$ 
and $\COL(\GI)$ depend on the subgraph $\GI$, whereas in Max-Cut the polytope $\CUT_{|I|}$ only 
depends on the number of vertices of $\GI$.

Finally let us mention that
an important feature of this hierarchy is that the size of the matrix variable remains $n$ or 
$n+1$ on all levels of the hierarchy.
On higher levels the ESCs are included into the SDPs in the most natural way
through convex combinations.
Hence on higher levels of the exact subgraph hierarchy new variables and linear 
constraints representing convex hull conditions are added to the SDP of the 
basic SDP relaxation.
 
Therefore it is 
possible to approximate $\eshmc{k}$, $\eshss{k}$ and $\eshc{k}$ by forcing only some subgraphs of 
order $k$ to be exact. This is our key ingredient to computationally obtain tight bounds on 
$z_{mc}$, $\alpha(G)$ and $\chi(G)$ and also a major advantage over several 
other SDP based 
hierarchies~\cite{LasserreHierarchy,LovaszSchrijverHierarchy,SheraliAdamsHierarchy}
 for NP-hard problems.

\subsection{Structural Differences of the Three Problems}
\label{sec:structuralDiff}
The focus of this paper lies in computational results, so we omit further extensive 
theoretical investigations, but we want to draw the attention to a major structural 
difference between the Max-Cut problem and the stable set and the coloring problem. 
Towards this end we consider  one graph from the Erd\H{o}s-R\'enyi 
model $G(n,p)$ with $n=100$ and $p = 0.15$. 
A graph from this model is 
a random graph of order $n$, in which each edge appears  
with probability $p$. 

We compute the optimal solutions of the basic 
relaxations~\eqref{relaxation mc},~\eqref{relaxation ss} 
and~\eqref{relaxation col} and denote them by $X^\ast$. Then for each subgraph 
$\GI$ of order 
$k\in \{2,3,4,5\}$ we compute the projection distance $\delta_{mc}^I$, 
$\delta_{ss}^I$ and $\delta_{c}^I$ of the 
submatrix 
$\XI^\ast$ of the corresponding $X^\ast$ to $\CUT_k$, 
$\STAB(\GI)$ and $\COL(\GI)$, respectively.
So  for example
\begin{align*}
\delta_{mc}^I = \min_{C \in \CUT_k} \norm{\XI^\ast - C},
\end{align*}
where $\norm{.}$ denotes the Frobenius norm. We consider a subgraph $\GI$ as 
violated, if the projection distance is larger than the small tolerance 
$5\cdot 10^{-5}$.

\begin{table}
    \setlength\tabcolsep{6pt}
    \centering
    \caption{Percentage of violated subgraphs of order $k$ for one random 
    graph}
    \begin{tabular}{|c|rrrr|}
        \hline
        $k$       & 2 &   3 &    4 &    5 \\ \hline
        \# subgraphs             & 4950& 161700& 3921225& 75287520 \\
        \% violated subgraphs MC & 0.00&  49.59&   91.69&    99.83 \\
        \% violated subgraphs SS & 7.54 &  21.96&   41.00&    60.88 \\
        \% violated subgraphs CO & 5.82&  16.83&   31.90&    49.14    \\       
        \hline
    \end{tabular}
    \label{tab:nrViolSubgraphs}
\end{table}

\begin{figure}
    \begin{center}
        \includegraphics[scale=0.08]{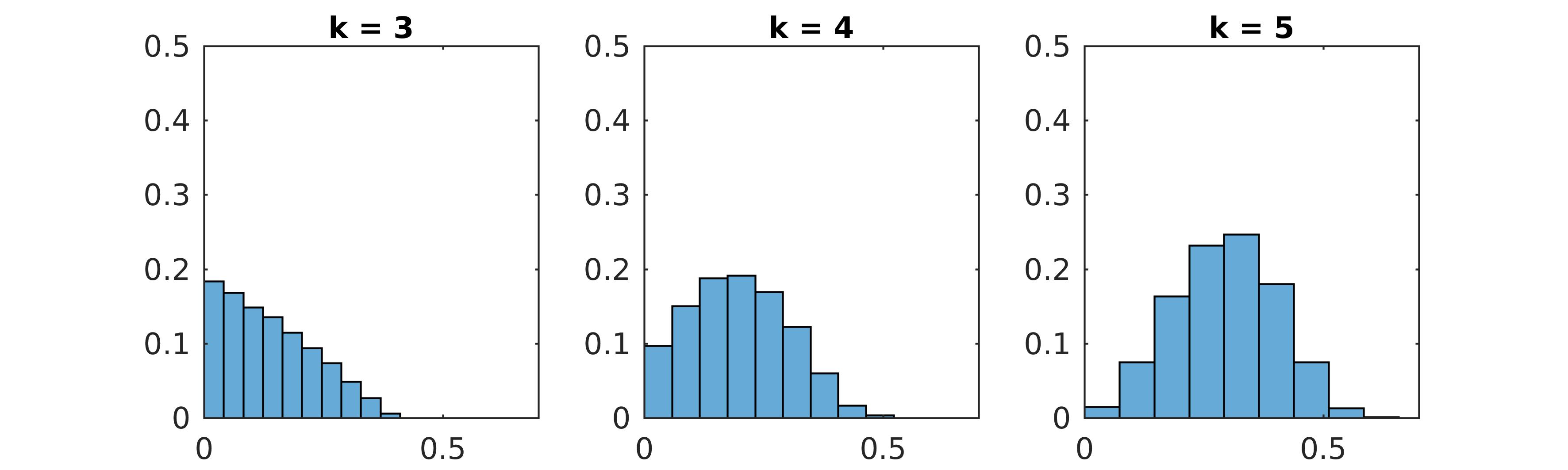}
        \caption{Histogram of $\delta_{mc}^I$ 
        for all violated 
            subgraphs $\GI$ of order $k\in \{3,4,5\}$}
        \label{fig:allProjDist mc}
    \end{center}
\end{figure}

\begin{figure}
    \begin{center}
        \includegraphics[scale=0.08]{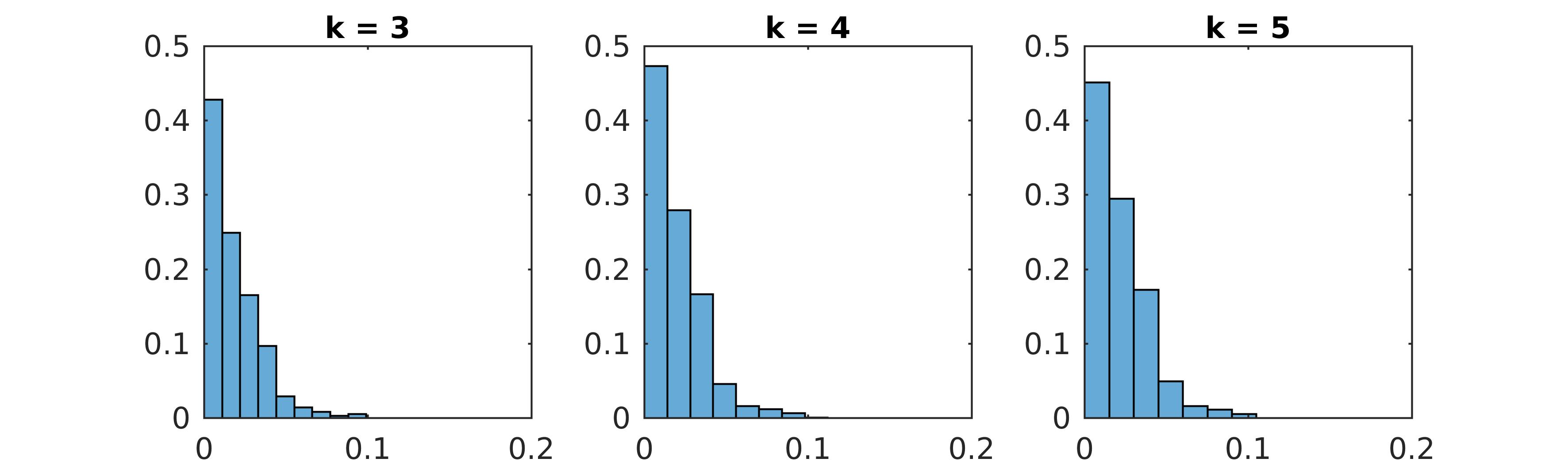}
        \caption{Histogram of $\delta_{ss}^I$
             for all violated 
            subgraphs $\GI$ of order $k\in \{3,4,5\}$}
        \label{fig:allProjDist ss}
    \end{center}
\end{figure}

\begin{figure}
    \begin{center}
        \includegraphics[scale=0.08]{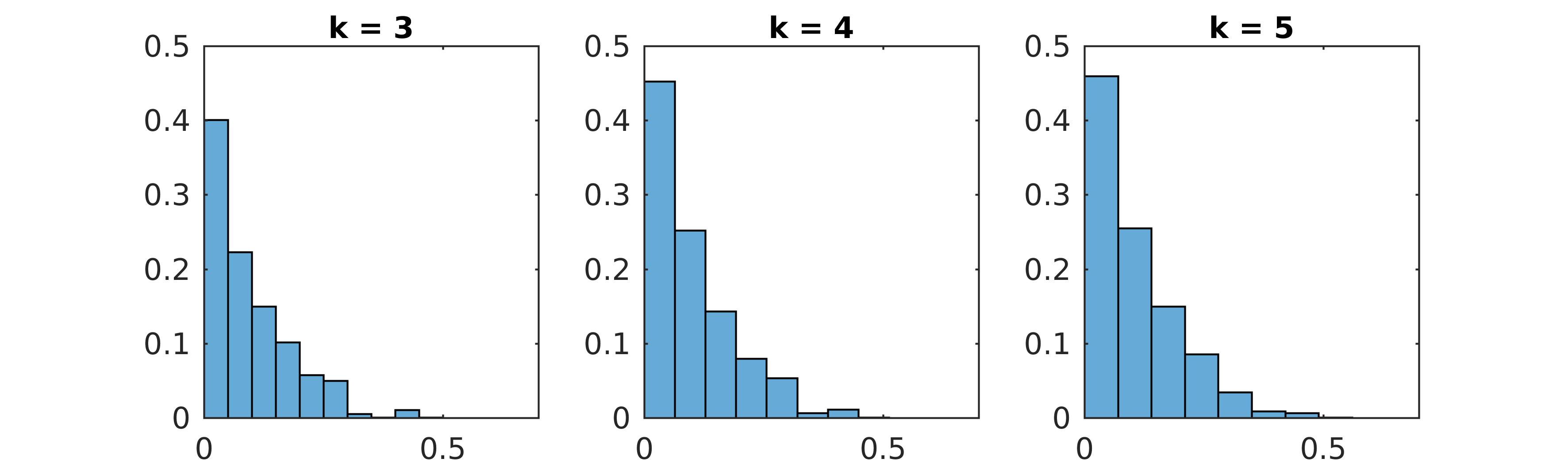}
        \caption{Histogram of $\delta_{c}^I$ for all violated 
            subgraphs $\GI$ of order $k\in \{3,4,5\}$}
        \label{fig:allProjDist co}
    \end{center}
\end{figure}

In Table~\ref{tab:nrViolSubgraphs} one sees that the number of violated 
subgraphs 
is much higher in the case of the Max-Cut problem than for the stable set and the 
coloring problem. Figure~\ref{fig:allProjDist mc},~\ref{fig:allProjDist ss} 
and~\ref{fig:allProjDist co} show the distribution of the 
projection distances 
of the violated subgraphs. 
They are normalized in such a way that $1$ is the total number of violated 
subgraphs.
Here it becomes obvious that for the Max-Cut problem most of the violated 
subgraphs have a large violation, whereas most of the violated subgraphs for 
the coloring problem have a small violation and an even smaller violation for 
the stable set problem.

Therefore in the case of the Max-Cut problem there are very many violated 
subgraphs, and typically all of them have a large projection distance. On the other hand 
for the stable set and the coloring problem only very few subgraphs have a 
large 
projection distance, the majority of the subgraphs is either not violated at all or only 
violated a little bit. Hence finding significantly violated subgraphs is much more 
difficult for the stable set and the coloring problem, than it is for the Max-Cut 
problem.

A possible explanation for this consists of the following
dimension argument.
Let $G$ be a graph on $n$ vertices with $m$ edges.
The SDP relaxation for Max-Cut
starts out with a matrix variable of size $n$ and $n$ equations, while the 
evaluation of $\vartheta(G)$ requires
a matrix of size $n+1$ and
$n+m+1$ equations
and in the computation of $t^{*}(G)$ there is a 
matrix of size $n+1$ and $2n + m$ equations. 
Hence the Max-Cut, stable set and coloring relaxation are contained in a 
$\binom{n}{2}$, $\binom{n}{2} + n - m$ and $\binom{n}{2} - m + 1$ dimensional 
space, and it makes sense that Max-Cut has the most and coloring has the least 
violated ESCs, just as we see it in Table~\ref{tab:nrViolSubgraphs}.
Furthermore in the stable set and the coloring relaxation the additional row 
and column together with the positive semidefiniteness constraint effect 
all
entries of $X$, even if they are not directly addressed by any constraint.
Therefore it is plausible that the violations for the Max-Cut problem are much 
larger than those for the stable set and the coloring problem.

For our computations that means that there is the hope that fewer ESCs are necessary to tighten the basic 
relaxation. This intuition is indeed confirmed in our computational experiments in Section~\ref{sec:comp}.


\section{Partial Lagrangian Dual} \label{sec: lagrangianDual}
We are interested in solving 
relaxations~\eqref{relaxation mc with esc},~\eqref{relaxation ss with esc} 
and~\eqref{relaxation col with esc} with a potentially large number of ESCs, where using interior 
point solvers is too time consuming. In this section we will first establish a unified 
formulation of the relaxations~\eqref{relaxation mc with esc},~\eqref{relaxation ss with esc} 
and~\eqref{relaxation col with esc}. Then we will build the partial Lagrangian dual of this 
formulation, where only the ESCs are dualized.

In order to unify the notation for the three problems observe that the ESCs $\XI \in \CUT_{|I|}$, 
$\XI \in \STAB^{2}(\GI)$ and 
$\XI \in \COL(\GI)$ can be represented as
\begin{equation}\label{esc in conv hull}
\XI = \sum_{i=1}^{\tI} \lambda_{i}C^{I}_{i},\quad \lambda \in \simplex_{\tI},
\end{equation}
where $C^{I}_{i}$ is the $i$-th cut, stable set or coloring matrix of the subgraph $\GI$ and 
$\tI$ is their total number.

A formal description of ESC in~\eqref{esc in conv hull} 
requires some additional notation. 
First we introduce the projection $\opP\colon \Sym{n} \mapsto \Sym{\kI}$, 
mapping $X$ to the submatrix $\XI$. 
Second we define a map $\opAmc\colon \Sym{\kI} \mapsto \R^{\tI}$, such that its adjoint map 
$\opAmct\colon \R^{\tI} \mapsto \Sym{\kI}$ is given by 
$\opAmct(\lambda)=\sum_{i=1}^{\tI}{\lambda_i C_i^{I}}$ and produces a linear combination of 
the cut, stable set or coloring matrices. 
Thus we can rewrite~\eqref{esc in conv hull} as 
\begin{equation}
\label{esc mc withAP}
\opAmct(\lambdaI) - \opP(X) = 0, \quad \lambdaI \in \simplex_{\tI}.
\end{equation}

The left-hand side of this matrix equality is a symmetric matrix, of which some entries (depending 
on which problem we consider) are zero for sure, so we do not have to include all 
$\kI \times \kI$ equality constraints into the SDP. Let $\bI$ be the number of equality 
constraints we have to include.
Note that $\bI = \binom{\kI}{2}$, $\bI = \binom{\kI}{2} + \kI - m_I$ and 
$\bI = \binom{\kI}{2} - m_I$ for the Max-Cut, stable set and coloring problem respectively, if
$m_I$ denotes the number of edges of $\GI$. This is because in the case of the stable set problem 
we also have to include equations for the entries of the main diagonal contrary to  Max-Cut and 
the coloring problem.
Then we define a linear map $\opMmc\colon \R^{\bI} \mapsto \Sym{\kI}$ such that the adjoint 
operator $\opMmct\colon \Sym{\kI} \mapsto \R^{\bI}$ extracts the $\bI$ positions, for which we 
have to include the equality constraints, 
into a vector. 
So we can rephrase~\eqref{esc mc withAP}
equivalently as
\begin{align*}
\opMmct(\opAmct(\lambdaI) - \opP(X)) = 0, \quad \lambdaI \in \simplex_{\tI},
\end{align*}
which are $\bI+1$ equalities and $\tI$ inequalities.
In consequence all three 
relaxations~\eqref{relaxation mc with esc},~\eqref{relaxation ss with esc} 
and~\eqref{relaxation col with esc}
have the generic form
\begin{equation}\label{sdp}
z = \max \{ \langle C, \widehat{X} \rangle:~
\widehat{X} \in \setX,~ \lambdaI \in \simplex_{\tI},~ \opMt(\opAt(\lambdaI) - \opP(X)) = 0  ~ \forall I \in J\},
\end{equation}
where $C$, $\setX$, $\opA$, $\opM$ and $\bI$ 
have to be defined in a problem specific way. Furthermore $\widehat{X} = X$ in 
the case of Max-Cut and 
stable set and $\widehat{X} = \left( 
\begin{array}{cc} t & \allones{}^{\opvt}\\ \allones{} & X 
\end{array}\right)$ for coloring, but for the sake of understandability we will just use $X$ in 
the following.

The key idea to get a handle on problem~\eqref{sdp} is to consider the partial
Lagrangian dual where the ESCs (without the constrains $\lambdaI \in \simplex_{\tI}$) 
are dualized. 
We introduce a vector of multipliers $\yI$ of size $\bI$ for each $I$ and collect them in 
$\yJ = (\yI)_{I\in J}$ and also collect $\lambdaJ = (\lambdaI)_{I\in J}$.
The  Lagrangian function becomes
\begin{equation*}
\Lagrangian(X,\lambdaJ,\yJ) = \langle C, X \rangle + \sum_{I \in J}{\langle \yI, \opMt(\opAt(\lambdaI) - \opP(X)) \rangle}
\end{equation*} 
and standard duality arguments (Rockafellar~\cite[Corollary 37.3.2]{RockafellarConvexAnalysis}) 
yield
\begin{equation}\label{lagrangianDual}
z = \min_{\yJ} \max_{\substack{X \in \setX\\ \lambdaI \in \simplex_{\tI}}} \Lagrangian(X,\lambdaJ,\yJ).
\end{equation}
For a fixed set of multipliers $\yJ$ the inner maximization becomes 
\begin{equation*} \max_{\substack{X \in \setX\\ \lambdaI \in \simplex_{\tI}}} \left\langle C - \sum_{I\in J}{\opPt \opM(\yI)}, X \right\rangle + \sum_{I\in J}{\langle \opA \opM(\yI), \lambdaI \rangle}.
\end{equation*}


This maximization is interesting in at least two aspects. 
First, it is separable in the sense that the first term depends only on $X$ and the second one 
only on the separate $\lambdaI$. 
Moreover, if we denote the linear map 
$\opA\opM\colon \R^{\bI} \mapsto \R^{\tI}$ 
with the matrix $\opD$, maximizing the summands of the second term is easy, 
because
the feasible region is a simplex. Hence the 
explicit solution of maximizing a summand of the second term is 
\begin{align}\label{solutionMaxTerm}
\max_{\lambdaI \in \simplex_{\tI}}\langle \opD (\yI), \lambdaI \rangle = 
\max_{1 \leqslant i \leqslant \tI}\left[\opD (\yI) \right]_{i}.
\end{align}

In order to consider the first term in more detail, we define the following function. Let 
$\sumbI = \sum_{I\in J}\bI$ be the dimension of $\yJ$. Then $h\colon \R^{\sumbI} \to \R$ is 
defined as
\begin{equation}\label{def:hy}
h(\yJ)= \max_{X \in \setX} \left\langle C - \sum_{I\in J}{\opPt \opM(\yI)}, X \right\rangle = \left\langle C - \sum_{I\in J}{\opPt \opM(\yI)}, X^{\ast} \right\rangle,
\end{equation}
where $X^{\ast}$ is a maximizer over the set $\setX$ for $y$ fixed. 
Note that $h(\yJ)$ is convex but non-smooth, but~\eqref{def:hy} shows that
\begin{equation}\label{def:nablah}
g_I= -\opM^T \opP(X^{\ast})
\end{equation}
 is a subgradient of $h$ with respect to $\yI$.  

With~\eqref{solutionMaxTerm} and~\eqref{def:hy} we reformulate the 
partial Lagrangian dual~\eqref{lagrangianDual} to
\begin{equation}\label{primal compact for bundle}
z = \min_{\yJ} \left\{ h(\yJ) + \sum_{I\in J}{\max_{1 \leqslant i \leqslant 
\tI}\left[\opD (\yI) \right]_i}\right\}.
\end{equation}

The dual formulation~\eqref{primal compact for bundle} of the original
semidefinite  
relaxation~\eqref{sdp}  has the form of a convex minimization problem
over the set of multipliers $y$. The evaluation of the function $h$ 
at a given $y$ requires solving a  ``simple'' SDP, independent of the 
number of ESCs included  in the relaxation.The function evaluation also
provides a subgradient of $h$ at $y$, given in \eqref{def:nablah}. 
Hence we propose to use the {\em bundle method} from convex optimization 
to solve~\eqref{primal compact for bundle}. The details are given in the 
subsequent section.


\section{Solving the Partial Lagrangian Dual}
\label{sec:bundle}

\subsection{The Bundle Method}

The bundle method is a well established tool in convex optimization to minimize a 
non-smooth convex function.
We refer to 
the recent monograph Bonnans, Gilbert,  Lemar{\'e}chal and 
Sagastiz\'{a}bal~\cite{NumericalOptimizationBonnans} for a nice introduction.
In our setting we want to use the bundle method in order to solve an SDP.  Helmberg and 
Rendl~\cite{HelmbergRendlBundleSDP} were the first to use a bundle method to solve SDPs 
in 2000. Later Fischer, Gruber, Rendl and 
Sotirov~\cite{FischerGruberRendlSotirovBundleMaxCutEquipartition} and Rendl and 
Sotirov~\cite{RendlSotirovBundleSDPQAP} used the bundle method for SDPs in order to get 
good relaxations for the Max-Cut and the equipartition problem and the 
quadratic assignment problem, respectively. 

The bundle method setting 
described by Frangioni and Gorgone 
in~\cite{EasyBundle}, which is set up to handle
$\max$ terms explicitly, is best suited for our purposes, so we
apply it to our problem~\eqref{primal compact for bundle}.

The bundle method is an iterative procedure. 
It maintains the \emph{current center}~$\yCenter$, 
representing the current estimate of the optimal solution, 
and the 
set $\Bundle = \{ 
(\yBdlExplanation{1}, \hBdlExplanation{1},\gBdlExplanation{1},\XBdlExplanation{1}), \dots, 
(\yBdlExplanation{\bdlSize}, 
\hBdlExplanation{\bdlSize},\gBdlExplanation{\bdlSize},\XBdlExplanation{\bdlSize})
 \}$, which is called  \emph{bundle},
 throughout the 
 iterations.
Here 
$y_1, \ldots, y_r$ are the points which we use to set up our subgradient 
model. Moreover $h_j = h(y_j)$, $g_j$ is a subgradient  of $h$ at $y_j$ 
and $X_j$ is a maximizer of $h$ at $y_j$ as in~\eqref{def:hy}. 

At the start we select $y_1=\yCenter =0$ and evaluate $h$ at $\yCenter$, which 
yields the bundle 
$\Bundle =\{(y_1, h_1,g_1,X_1)\}$. A general iteration consists 
of first determining the 
new \emph{trial point}, then
evaluating the function at this new point, and finally updating the bundle
$\Bundle$.
In the literature evaluating the function is referred to as calling the 
\emph{oracle}.
The subgradient information of the bundle $\Bundle$ 
translates into the subgradient model
\begin{align*}
h(y) \geqslant \hBdlExplanation{j} + \langle \gBdlExplanation{j},y-\yBdlExplanation{j}\rangle \text{ for all } j = 1, \dots,  \bdlSize.
\end{align*}
It is common to introduce
$$
\eBdlExplanation{j} = h(\yCenter) - \hBdlExplanation{j} - \langle \gBdlExplanation{j},\yCenter-\yBdlExplanation{j}\rangle \text{ for } j = 1, \dots,  \bdlSize
$$
and to define $e = (e_j)_{j=1,\dots,r}$.
With $\hCenter = h(\yCenter)$ the subgradient model becomes
\begin{align}
\label{subgrMaxModel}
h(y) \geqslant \max_{1\leqslant j \leqslant r} \left\{\hCenter - \eBdlExplanation{j} + \langle \gBdlExplanation{j},y-\yCenter\rangle \right\}.
\end{align}
The right-hand side above is convex, piecewise linear and minorizes $h$. 
In each iteration of the bundle method we minimize the right-hand 
side of~\eqref{subgrMaxModel} instead of $h$, 
but ensure that we do not move too far from $\yCenter$ 
by adding a penalty term 
of the form
$
\frac{1}{2}\bdlPar\norm{y-\yCenter}_2^{2}
$
for a parameter $\bdlPar \in \R_{+}$  
to the objective function.
We introduce auxiliary variables $w \in \R$ and $\vI \in \R$ for all $I \in J$ 
to model the maximum 
terms. With  $\nrESC = |J|$ and $v = (\vI)_{I \in J} \in \R^\nrESC$ we  
end up with
\begin{align}
\label{bundleProblemPiO2}
\min_{y,w,v} \quad w &+ \sum_{I \in J} \vI + \frac{1}{2}\bdlPar\norm{y-\yCenter}_2^{2}\\ \nonumber
\st \quad w &\geqslant \hCenter - \eBdlExplanation{j} + \langle 
\gBdlExplanation{j},y-\yCenter\rangle && \forall j = 1,\dots,\bdlSize \\
\nonumber
\vI &\geqslant \left[\opD (\yI) \right]_{i} && \forall i = 1,\dots,\tI \quad 
\forall I \in J.
\end{align}
This is a convex quadratic problem in $1+\nrESC+\sumbI$ variables 
with $\bdlSize+\sum_{I \in J}\tI$ linear inequality constraints
which is often referred to as the \emph{bundle master problem}. 
Its solution $(\yTrial,\wTrial,\vTrial)$ provides 
the new trial point $\yTrial$. In the following section we will
discuss computational issues and present a practically 
 efficient approach
starting with its dual, see below.  

The second step in each bundle iteration is to evaluate the 
function $h$ at $\yTrial$ which means solving 
the basic SDP relaxation as introduced in 
Section~\ref{sec:DefProblemsBasicRel} with a modified objective function. 
In the case of 
Max-Cut this function evaluation 
can be done very 
quickly (solve an SDP with $n$ simple equations).  
For the stable set and the coloring problem 
the resulting SDP is computationally
more demanding, as there are also equations for each edge in the graph. 
The bundle iteration is finished by deciding whether $\yTrial$ becomes the 
    new center (serious step, roughly speaking if the increase of the objective 
    function is good enough) 
    or not (null step). 
    In either case the new point is included in the bundle, 
    some other elements of the bundle are possibly removed, the bundle parameter $\bdlPar$ is 
    updated  
    and a new iteration starts.
    
\subsection{The Dual of the Bundle Master Problem}

In the bundle method it is commonly proposed to solve the dual problem 
of~\eqref{bundleProblemPiO2}, hence next we derive the dual of~\eqref{bundleProblemPiO2}.
Towards this end we collect the subgradients $\gBdlExplanation{i}$ in the 
matrix $\opG=(\gBdlExplanation{1}, \ldots, \gBdlExplanation{\bdlSize})$. 
It will be notationally convenient to partition the matrix $\opG$ into blocks of rows 
corresponding to the subsets $I \in J$, namely 
$\opG = (\opGI)_{ I \in J }$ where each $\opGI$ has $r$ columns and $\bI$ 
rows.
Furthermore we make the subgradient model and maximum term constraints more compact by 
reformulating them to 
$w \allones{} \geqslant \hCenter\allones{} - e + \sum_{I \in J} \opGItrans(\yI -\yICenter)$ 
and $\vI\allones{} \geqslant \opD(\yI)$. 

We denote by $\alpha \in \R^{\bdlSize}$ the dual variables to the subgradient model 
constraints and with $\betaI \in \R^{\tI}$ the dual variables of the constraints 
involving $\vI$ for the maximum terms. Furthermore we define 
$\beta = (\betaI)_{I \in J}$ as the collection of all $\betaI$. 
Hence we obtain the Lagrangian function
\begin{align*}
\Lagrangian(y,w,v,\alpha,\beta) = w &+ \sum_{I \in J}\vI + 
\frac{1}{2}\bdlPar \sum_{I \in J}\norm{\yI-\yICenter}_{2}^{2} \\ 
& \quad \quad + \left\langle \alpha, \hCenter \allones{} -e - w\allones{} \right\rangle +
\sum_{I \in J} \left\langle \alpha, \opGItrans (\yI -\yICenter) \right\rangle \\
& \quad \quad + \sum_{I \in J} \left\langle \betaI, \opD( \yI) - \vI \allones{} 
\right\rangle. 
\end{align*}
After exchanging $\min$ and $\max$ by using 
strong duality the dual of~\eqref{bundleProblemPiO2} becomes
\begin{align*}
\max_{\substack{\alpha\geqslant 0\\\beta \geqslant 0}} \text{ } \min_{y,w,v} \Lagrangian(y,w,v,\alpha,\beta).
\end{align*}
Since $\nabla_{w}\Lagrangian=0$,~ 
$\nabla_{\vI}\Lagrangian = 0$, 
and~$\nabla_{\yI} \Lagrangian = 0$ has to 
hold for all $I \in J$ at the dual optimum, we get $\alpha \in \simplex_{\bdlSize}$, 
$\betaI \in \simplex_{\tI}$ and 
\begin{align}
\label{formulaOptimaly}
\yI = \yICenter - \frac{1}{\bdlPar}\left(\opGI(\alpha) + \opDt(\betaI)\right).
\end{align} 
In consequence the dual of~\eqref{bundleProblemPiO2} simplifies to
\begin{align}
\label{dPiO2}
\max_{\substack{\alpha\in \simplex_{\bdlSize}\\ \betaI \in \simplex_{\tI}}}
\hCenter - e^{\transposedVec} \alpha + \sum_{ I \in J }\left\langle \opD 
(\yICenter) , \betaI \right\rangle - 
\frac{1}{2\bdlPar}\sum_{I \in J} \norm{\opGI(\alpha) + \opDt(\betaI)}^{2}_{2}.
\end{align}
This is a convex quadratic problem with $r + \sum_{I\in J} \tI$ variables 
and $1+\nrESC$ simple equality constraints, asking that the respective block of
variables adds up to one.
Now instead of solving~\eqref{bundleProblemPiO2} within the bundle method 
directly, we solve its dual~\eqref{dPiO2} to 
get the multipliers $\alpha$ and 
$\beta$ and recover $\yTrial$ using~\eqref{formulaOptimaly}.

\subsection{Our Bundle Method}

So far we have sketched how to use our bundle method in order 
to obtain a solution $y$ of~\eqref{primal compact for bundle}, 
but actually 
we are interested in a solution $X$ of~\eqref{sdp}. 
One can use the bundle 
$\Bundle = \{ 
(\yBdlExplanation{1}, 
\hBdlExplanation{1},\gBdlExplanation{1},\XBdlExplanation{1}), \dots, 
(\yBdlExplanation{\bdlSize}, 
\hBdlExplanation{\bdlSize},\gBdlExplanation{\bdlSize},\XBdlExplanation{\bdlSize})
\}$,
which is updated in each iteration, 
in order to obtain a good 
approximate solution for $X$. In particular 
it follows from the convergence 
theory of the bundle method that under mild conditions
\begin{align}\label{eq:optXlam}
X = \sum_{j=1}^{r} \alpha_j X_j \quad \text{and}\quad  \lambdaI = \betaI
\end{align}
converges to the optimal values of $X$ and $\lambdaI$ of~\eqref{sdp}, 
see for example Robinson~\cite{RobinsonBundlePrimalConvergence} for the general 
theory and Gaar~\cite{elli-diss} for the convergence in our particular 
setting.

We are now able to present our version of the bundle method.
Note that there is no need of keeping $y_j$ in the bundle explicitly 
by computing and updating $e$ in a proper way, 
so we drop $y_j$ from the bundle $\Bundle$.
Algorithm~\ref{alg: general procedurce bundle} summarizes the main 
computational steps 
of our bundle method  to get an approximate
optimal solutions of~\eqref{sdp} and~\eqref{primal compact for bundle}.

\begin{algorithm}[ht]
    \caption{Solving \eqref{sdp} and \eqref{primal compact for bundle} for a 
    given $J$}
    \label{alg: general procedurce bundle}
    \begin{algorithmic}[1]
        \Procedure{Our\_Bundle}{$h$, $\opD$}
        \State $\overline{y}=y_1 = 0$
        \State Evaluate $h$ at $y_1$ to get
        $h_1 = h(y_1)$, 
        a maximizer $X_1$ and a subgradient $g_1$
        \State Set the bundle to $\Bundle = \{(h_1, g_1, X_1)\}$
        \While{``stopping condition is not satisfied''}
        \State Solve the convex quadratic problem \eqref{dPiO2} to get 
        $\alpha$ and $\beta = (\betaI)_{I \in J}$
        \State Determine $\yTrial$ using~\eqref{formulaOptimaly}
        \State Determine $X$ and $\lambda$ using~\eqref{eq:optXlam}
        \State Evaluate $h$ at $\yTrial$ to get $h(\yTrial)$, a maximizer 
        $\widetilde{X}$ and a subgradient $\widetilde{g}$
        \State Decide whether $\yTrial$  becomes the new center $\overline{y}$ 
        (serious step) 
        or not (null step)
        \State Update the bundle $\Bundle$ and the bundle parameter $\mu$
        \EndWhile
        \State \Return $\yCenter$, $X$, $\lambda$
        \EndProcedure
    \end{algorithmic}
\end{algorithm}

The generic description of our bundle method in 
Algorithm~\ref{alg: general procedurce bundle} leaves some flexibility to the 
user. 
We will present implementation details in 
Section~\ref{sec:bundleImplementationDetails}.


\section{The Computation of ESCs Based Bounds}
\label{sec:overallAlgo}
\subsection{The Overall Algorithm}
The goal of this paper is to get good bounds on the optimal Max-Cut 
value $z_{mc}$, the stability number $\alpha(G)$ and the chromatic number 
$\chi(G)$ by including ESCs into the basic SDP 
relaxations~\eqref{relaxation mc},~\eqref{relaxation ss} 
and~\eqref{relaxation col} 
in order to improve the bounds from the basic SDP relaxations. 
We will call bounds obtained in this way \emph{exact subgraph bounds (ESB)}.
In other words ESBs are attained by 
solving~\eqref{relaxation mc with esc},~\eqref{relaxation ss with esc} 
and~\eqref{relaxation col with esc} or, in the generic form, by 
solving~\eqref{sdp}.

Up to now we have concentrated on the most subtle part of retrieving good ESBs, 
which consists in solving the SDP relaxation~\eqref{sdp} with a given set $J$ 
of ESCs. 
Our ultimate goal however is to reach ESBs where all ESCs
of order $k$ are (nearly) satisfied for small values of $k$ like 
$k\leqslant 7$.

We propose to reach this goal by proceeding iteratively. 
Starting with $k=3$ in the Max-Cut case 
(as there are no violated ESCs of order 2) and $k=2$ in the other cases
we search for violated ESCs of order $k$ and include only the most
violated ESCs that we find into $J$. 
After solving the SDP~\eqref{sdp}, 
we follow an extreme strategy and remove any ESC that has become
inactive. As we typically still find further badly violated ESCs this 
allows us a quick exploration of the entire space of ESCs. 
Once we do not find ESCs of order $k$ with significant violation, we
increase $k$ and continue. 
We call each such iteration a
\emph{cycle}. 

In each cycle so we keep some information, such as 
the current dual variables $y_i$ and the bundle $\Bundle$,  appropriately 
modified to reflect possibly deleted and added new constraints. 
In particular we delete from all $y_i$ the positions corresponding to deleted 
ESCs, extend all $y_i$ with zeros for the newly added ESCs and 
deduce the update of all other variables. This choice allows us to reuse the 
bundle $\Bundle$. 
Our procedure to compute ESBs is sketched in 
Algorithm~\ref{alg: general procedurce bound computation}.

\begin{algorithm}[ht]
    \caption{Computation of an exact subgraph bound}
    \label{alg: general procedurce bound computation}
    \begin{algorithmic}[1]
        \Procedure{Compute\_ESB}{$G$}
        \State $k = 2$ (or $k=3$ for Max-Cut)
        \State $J = \emptyset$
        \While{``stopping condition is not satisfied''}
        \State Get an approximate solution $X$ of~\eqref{sdp} and $y$  
        of~\eqref{primal compact for bundle} \label{alg:lineApproxSol}
        \label{algStep:bdlMasterProblem}       
        \State Update the ESB to the objective function value of $y$
        of~\eqref{primal compact for bundle}          
        \State Remove inactive ESCs from $J$
        \State Include most violated ESCs of $X$ with order $k$ into $J$ 
        \If{``not enough violated subgraphs found''}
        \State $k = k + 1$
        \EndIf      
        \EndWhile
        \State \Return ESB
        \EndProcedure
    \end{algorithmic}
\end{algorithm}

\begin{figure}
    \begin{center}
        \includegraphics[width=0.5\linewidth]{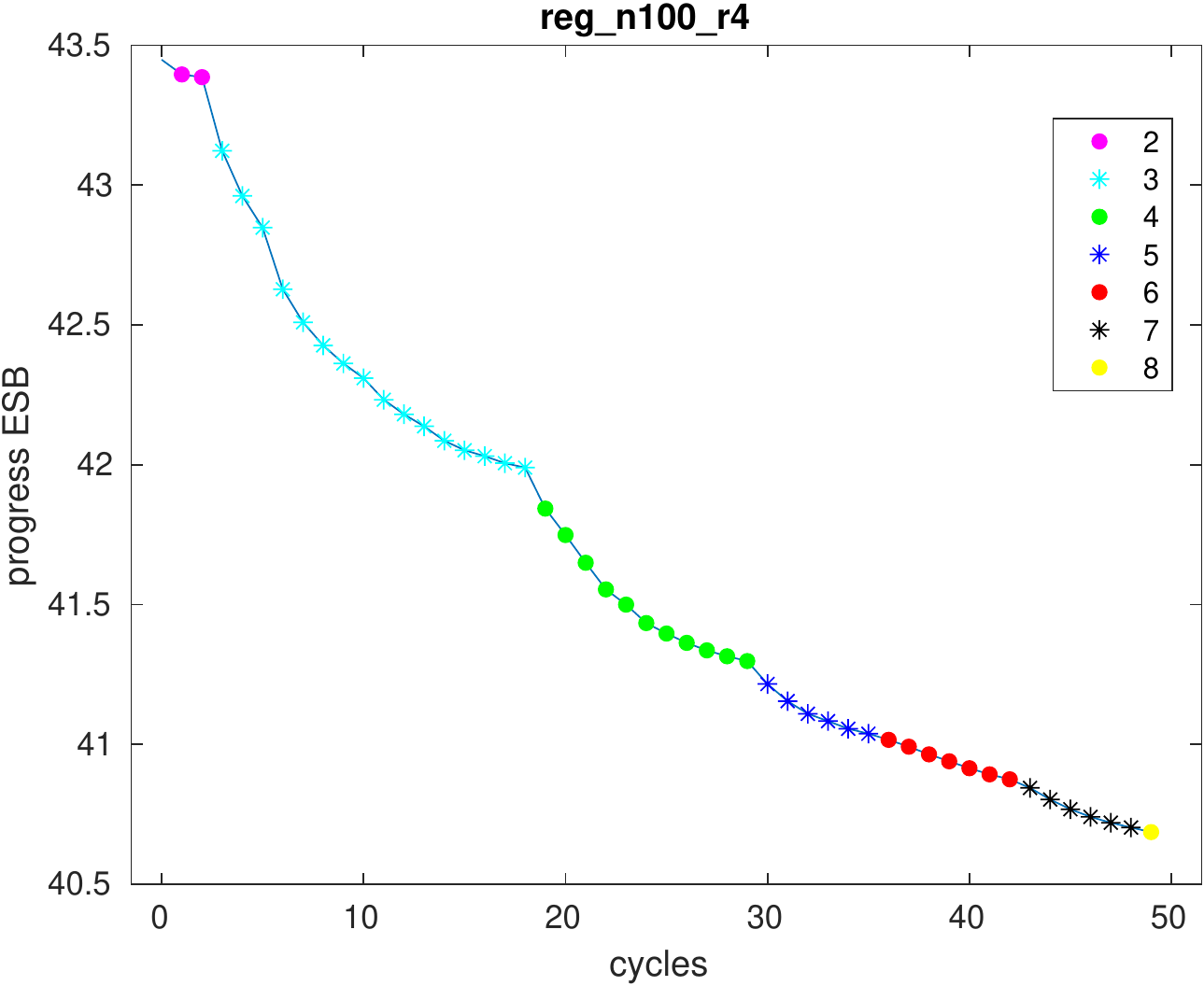}
        \caption{Progress of the ESB over 50 cycles for one instance of 
        Table~\ref{regulargraphs}}
        \label{fig:progressttf}
    \end{center}
\end{figure}

The typical behavior over a set of cycles for one stable set instance can be 
seen in 
Figure~\ref{fig:progressttf}. After only a few cycles with $k=2$ we 
move to $k=3$. Here it takes 16 cycles to reach a point with 
all ESCs nearly satisfied. The Figure clearly shows a continuing
improvement of the ESB over the cycles. 

Note that the ESB computed in Algorithm~\ref{alg: general procedurce bound 
computation} is indeed a valid bound, because any $y$ is feasible 
for~\eqref{primal 
compact for bundle} and hence its dual objective function value is a valid 
bound on the primal optimal objective function value~\eqref{sdp}, which in turn 
is a valid bound on the optimal objective function value of the combinatorial 
optimization problem. Hence it is not necessary 
to solve~\eqref{sdp} 
and~\eqref{primal compact for bundle} to optimality to obtain valid bounds.
Of course we want to use our 
bundle method, Algorithm~\ref{alg: general procedurce bundle}, in order to 
obtain the approximate solutions in 
line~\ref{alg:lineApproxSol} of 
Algorithm~\ref{alg: general procedurce bound computation}.

\subsection{Finding Violated Exact Subgraph Constraints}
The key ingredients of 
Algorithm~\ref{alg: general procedurce bound computation} are on the one hand 
Algorithm~\ref{alg: general procedurce bundle}, which was detailed in 
Section~\ref{sec:bundle}, and on the other hand the update of the set $J$. 

The crucial point in order to do so is to find violated ESCs.
Let $\GI$ be a subgraph of oder $\kI$ of $G$ and 
$X^\ast$ be the current solution 
of~\eqref{sdp} and let $U$ be an arbitrary  
$\kI \times \kI$ matrix. 
Clearly $\CUT_{\kI}$, $\STAB^2(\GI)$ and $\COL(\GI)$ are bounded 
polytopes, hence the inner product of any element of these polytopes 
with $U$ is contained in a certain interval. 
Thus finding $I$ such that the inner product of $U$ with the submatrix 
$\XI^\ast$ of $X^\ast$ is minimum identifies  
a potentially violated subgraph. 

This minimization may be recast as a quadratic assignment 
problem consisting of the data matrices $X^{*}$ and the matrix $U$
embedded in an $n \times n$ matrix.  
We repeatedly 
use a local search heuristic for different fixed $U$ in order to obtain 
potentially violated subgraphs.
Then we compute the projection distances of $\XI^\ast$ to the corresponding 
polytope for all these subgraphs $\GI$ and include those into $J$ which have 
the largest projection distances and hence are violated most.

Possible choices for $U$ make use of 
hyperplanes for the respective target polytope, 
but other choices are possible. 
In our computations we use a collection of different matrices for $U$, for 
example matrices that induce facets of the corresponding polytope (if their 
computation for a particular $\kI$ is possible easily, which is the case for 
$\kI \leqslant 6$), extreme copositive matrices 
with $\{0,1,-1\}$ entries and random matrices. For each cycle we use at most 
$50$ different matrices $U$.

\subsection{Details of the Bundle Implementation}
\label{sec:bundleImplementationDetails}

We now briefly discuss some details of our implementation of 
Algorithm~\ref{alg: general procedurce bundle}
when used in 
line~\ref{alg:lineApproxSol} of 
Algorithm~\ref{alg: general procedurce bound computation}.
First of all one needs to decide on a stopping condition. 
Ideally we would stop, once a subgradient equal to zero is found. 
In our case, we either stop once the norm of the new subgradient is small 
enough (in the case of Max-Cut), or once the difference of the value of the 
function at the current center point and the value of the subgradient model of 
the 
function at 
the new trial point is smaller than some tolerance (as it is done 
in~\cite{NumericalOptimizationBonnans}, the tolerance is 0.005 in our 
implementations) 
or once we reach a maximum number of iterations (30 in our implementations). 
The third condition 
is motivated by the fact that we typically will continue adding new violated
ESCs, so there is no real need to get the exact minimum 
of~\eqref{primal compact for bundle}. 
Note however that it is important to come  
close to the optimal solution, because otherwise the 
resulting $X$ does not have a 
high 
enough precision in order to be useful for finding new violated subgraphs.

For updating the bundle we always add the new trial point to the bundle, 
but remove subgradients from the bundle that have become inactive.
This extreme choice of updating the bundle led to the best performance in our 
computational experiments.
In order to update the bundle parameter $\mu$ we use a modification of an 
update proposed by Kiwiel~\cite{KiwielProximityControl}.
We perform a serious step whenever the improvement of the objective function 
value of the new trial point is at leas a certain fraction of the expected 
improvement. This is a standard criterion, see for example~\cite{HUL2}. 
We solve the bundle master problem as a rotated second-order cone program 
(see~\cite{IntroductionSOCP} for more details) with MOSEK.


\section{Computational Results}
\label{sec:comp}

\subsection{Bundle Approach versus Interior Point Methods}

We start our computational investigation with a comparison of our bundle method 
with an interior point method in order to solve~\eqref{sdp}. 
In our overall Algorithm~\ref{alg: general procedurce bound computation} 
presented in Section~\ref{sec:overallAlgo} this has to 
be done in each cycle, so we are highly interested in fast running times.

From a theoretical point of view it is clear which method will win this 
competition: 
Assume we include $q = 1000$ ESCs (so $q = |J|$) for subgraphs of order $\kI = 
5$ 
in~\eqref{sdp} for the stable set problem. Then we have $\tI \leqslant 2^5 = 
32$ stable set matrices that potentially span $\STAB^2(\GI)$, and up to $\bI 
\leqslant \binom{\kI}{2} + \kI = 15$ equality and one inequality 
constraint for each ESC. In total we have up to $32000$ variables 
that have to fulfill up to $16000$ constraints in~\eqref{sdp}---additionally 
to 
the variables and constraints of the basic SDP 
relaxation~\eqref{relaxation ss}.
It is clear that the number of constraints will be a challenge for 
an interior point solver. 
In particular an interior point solver has to solve this 
SDP with a large number of constraints at once, whereas our 
bundle method in Algorithm~\ref{alg: general procedurce bundle} ``only''
has to solve the basic SDP relaxation and the bundle master problem over 
several iterations.
Therefore, we expect the bundle method to be the clear winner in this 
competition and refrain from a large scale comparison.

Instead, we compare the two methods only on some instances to confirm our 
theoretical inspection.
In Table~\ref{MC and SS solution times exact} we list the results for one 
Max-Cut and one stable set instance, both are taken from the 
Erd\H{o}s-R\'enyi model $G(n,p)$. 
We vary the number of included ESCs for subgraphs of order $3$, $4$ and $5$, so 
we solve~\eqref{relaxation mc with esc} and~\eqref{relaxation ss with esc} 
for different $J$. 
We choose $J$ such that the total number of equality constraints induced 
by the convex hull formulation of the ESCs $b$ ranges between 6000 and 15000. 
On the one hand we solve the instances with two interior point solvers, namely 
MOSEK and SDPT3~\cite{SDPT31,SDPT32} and list the running times in seconds.
On the other hand we use our bundle method. In our context we are mostly 
interested to improve
the upper bounds quickly, so we do not run 
Algorithm~\ref{alg: general procedurce bundle} until
we reach a minimizer, but stop after $30$ iterations. We list the running time 
for the oracle, i.e.\ the sum of the solution times of the basic SDP 
relaxation, and the overall running times. Additionally we present how much 
$\%$ of the MOSEK running time the bundle method needs and how close the 
solution found by the bundle method is to solution of MOSEK in $\%$ ($100\%$ 
means the solutions coincide).

In Table~\ref{MC and SS solution times exact} one sees that the running times 
decrease drastically if we use the bundle method compared to interior point 
solvers. For $b\approx 15000$ it takes the bundle method only around $8\%$ of 
the MOSEK running time to get as close as $95\%$ to the optimal value, which is 
sufficient for our purposes. One sees that our bundle method scales much better 
for increasing $|J|$, so for an increasing number of ESCs.
Furthermore MATLAB requires 12 Gigabyte of memory with interior point solvers 
for $b=15000$, 
showing also memory limitations. 

\begin{table}
    \setlength\tabcolsep{4.5pt}
    \centering
    \caption{Running times for one Max-Cut and one stable set instance with 
    different  
        sets of ESCs, where the graphs of order $n=100$ are from the 
        Erd\H{o}s-R\'enyi model}
    \begin{tabular}{|c|rrr|r|rr|rr|rr|}
        \hline
        & & & & & \multicolumn{2}{c|}{interior point} & \multicolumn{4}{c|}{
            Algorithm~\ref{alg: general procedurce bundle}}
        \\ \hline        
        & \multicolumn{3}{|c|}{\#ESC of order} & &\multicolumn{2}{c|}{time 
        (sec)}
        & \multicolumn{2}{c|}{time (sec)}   & \multicolumn{2}{c|}{$\%$ of MOSEK}
        \\\hline
        & $3$  & $4$  &  $5$   &  $b$        & MOSEK   & SDPT3     
        & oracle & overall & time & value\\\hline
        \multirow{ 4}{*}{MC}& 2000  & &    0 &  6000 &    18.37 &   49.22 &  
        1.01 &  6.05 & 32.93 & 97.20\\
        & 2000  & &  300 &  9000 &    55.24 &  134.78 &  1.18 &  9.33 & 16.90 & 
        95.02\\
        & 4000  & &    0 & 12000 &   104.56 &  289.78 &  1.71 & 11.13 & 10.64 & 
        93.66\\
        & 3000  & &  600 & 15000 &   184.43 &  525.85 &  1.56 & 14.83 & 8.04  & 
        94.54\\        
        \hline 
        \multirow{ 4}{*}{SS}& 1050 &    0 &    0 &  5914 &   23.54 &   79.25 &  
        7.86 & 10.65 & 45.22 & 98.25 \\
        & 1050 &  212 &   63 &  8719 &   50.11 &  174.33 & 10.61 & 16.52 & 
        32.96 & 97.89\\
        & 2100 &    0 &    0 & 11780 &  126.40 &  388.07 &  7.43 & 12.27 & 9.71 
        & 93.65 \\
        & 1575 &  318 &  212 & 14653 &  241.29 &  648.83 & 10.79 & 20.21 & 8.38 
        & 94.44\\
        \hline
    \end{tabular}
    \label{MC and SS solution times exact}
\end{table}

To summarize our small computational investigation confirms our intuition that 
the bundle method is much better suited for our purposes.

We want to point out that the number of bundle iterations can be increased in 
order to get closer to the optimum. For the 
larger instances in Table~\ref{MC and SS solution times exact} this will still 
result in significantly shorter running times.

Note that the bundle method has another advantage: A warm start with the bundle 
$\Bundle$ and the solution $\yCenter$ of the previous iteration in 
line~\ref{algStep:bdlMasterProblem} of Algorithm~\ref{alg: general procedurce 
bound computation} is possible. Since many ESCs remain the same in $J$  the 
problem to solve in 
line~\ref{algStep:bdlMasterProblem} does not change too much and a warm start 
can be very beneficial.

As a last remark we want to draw the attention to the running times for the 
oracle in Table~\ref{MC and SS solution times exact}. For the stable set 
problem the oracle needs over half of the running time, whereas in the Max-Cut 
problem the oracle evaluation is much faster. This is due to the fact that the 
basic SDP relaxation is a simpler SDP for the Max-Cut problem. 

In the following we present several computational results for obtained ESB by 
using the bundle method. Note that we refrain 
from comparing the running times of our bundle method with the running 
time of interior point 
methods, because interior point methods would reach their limit very soon.

\subsection{The Stable Set and the Coloring Problem}

In this section we will extend the computational results from \cite{GaarRendl} 
for the stable set and the coloring problem.
The computational investigations show that 
(i) the ESB obtained by 
including ESCs of fixed order $\kI$ improve for increasing $\kI$ and 
(ii) after including 
several ESCs for subgraphs of order $\kI$ the maximum projection distance of 
the 
violated subgraphs found decreases drastically. 

We extend these computational results by deriving one final ESB for 
several instances with  
Algorithm~\ref{alg: general procedurce bound computation}. We stop as soon as 
we have performed 50 cycles and only include subgraphs of order 
$k \leqslant 8$. We add at most $100$ ESCs 
in each cycle and warmstart the bundle with the information of the previous 
cycle. We already saw the typical behavior of the ESB over the cycles in 
Figure~\ref{fig:progressttf}.

\begin{table}
    \begin{center}
        \setlength\tabcolsep{5pt}
        \caption{Stable set results for torus graphs}
        \label{torusgraph}
        \begin{tabular}{|rrr|rrr|rr|} \hline
            $ d$ & $n$ & $m$ & $\vartheta(G)$ 
            & \shortstack{final \\ ESB} & $\alpha(G)$ 
            & \shortstack{time\\oracle} & \shortstack{time \\ other} \\ \hline
5 &  25 &     50 &  11.180 &  10.002 &  10  &       14.45 &   18.58\\ 
7 &  49 &     98 &  23.224 &  21.009 &  21  &       95.06 &   31.96 \\ 
9 &  81 &    162 &  39.241 &  36.021 &  36  &      277.98 &   66.37 \\ 
11 & 121 &    242 &  59.249 &  55.066 &  55 &      728.78 &  122.46 \\ 
13 & 169 &    338 &  83.254 &  79.084 &  78 &      859.73 &  170.85 \\ 
15 & 225 &    450 & 111.257 & 106.287 & 105 &     1390.52 &  224.97 \\ 
17 & 289 &    578 & 143.259 & 136.821 & 135 &     3123.41 &  314.29 \\  
\hline
        \end{tabular}
    \end{center}
\end{table}

As a first structural easy class of graphs we consider two-dimensional torus 
graphs which
are constructed as follows. For given  $d$, the graph $T_d$
has 
$d^2$ vertices which we label by $(i,j)$ for $i,j \in \{1, \ldots,d\}$.
The vertical edges join vertices with neighboring $i$ indices (and $j$ fixed),
yielding  edges $\{(i,j), (i+1,j)\}$ modulo $d$, and similarly the horizontal 
edges join vertices with $i$ fixed $\{(i,j), (i,j+1)\} $ modulo $d$. So there 
is a total
of $n=d^2$ vertices and $m=2n$ edges. It is not hard to verify that
in case of odd $d=2t+1$, we get $\alpha(T_d)=t(2t+1)$ and if $d=2t$ we have
$\alpha(T_d) = 2t^2$. The even case is not interesting, as $\vartheta(T_d)= 
\alpha(T_d)$. For $d$ odd we summarize some computational results in  
Table~\ref{torusgraph}.
We observe that for these graphs our ESB is substantially better than 
$\vartheta(G)$ and we close the integer gap for all  
instances with $n \leqslant 121$.

When considering the running times observe that the majority of the running 
time (given in seconds in Table~\ref{torusgraph}) is used for the oracle, 
because the SDP to evaluate $\vartheta(G)$ given in~\eqref{relaxation ss} 
with a 
slightly modified objective function is nontrivial. 
We tried several solver to 
solve this SDP, among them the interior point solver MOSEK~\cite{mosek}, and 
solvers based on alternating direction method of multipliers as 
DADAL~\cite{DADAL} and SDPNAL+~\cite{SDPNALPlus}. Both these solvers show very 
good results on computing $\vartheta(G)$, but as soon as the objective function 
slightly changes they do not perform well anymore.
Hence it will be future research to develop an SDP solver dedicated to these 
kind of instances.
Note that the running time in order to perform Algorithm~\ref{alg: general 
procedurce bound computation} is not very high and in particular 
only  increases mildly for larger instances.

\begin{table}
    \begin{center}
        \caption{Stable set results for near-regular graphs}
        \label{regulargraphs}
        \setlength\tabcolsep{5pt}
        \begin{tabular}{|rrrr|rrc|rr|} \hline
            graph  & $n$ & $r$ & $m$ & $\vartheta(G)$ 
            & \shortstack{final\\ ESB} & $\alpha(G)\geqslant$ 
            & \shortstack{time\\oracle} & \shortstack{time \\ other} \\ \hline
reg\_n100\_r4 & 100 & 4& 195 &  43.449 &  40.687 &  40 & 1020.60 &  143.12\\ 
reg\_n100\_r6 & 100 & 6& 294 &  37.815 &  35.246 &  34 &  935.69 &  125.93\\ 
reg\_n100\_r8 & 100 & 8& 377 &  34.480 &  32.190 &  31 &  939.27 &  127.46\\ 
reg\_n100\_r10 & 100&10& 474 &  31.797 &  28.954 &  28 & 1126.75 &  140.42\\ 
\hline 
reg\_n200\_r4 & 200 & 4& 400 &  87.759 &  83.732 &  80 & 1278.61 &  158.67  \\ 
reg\_n200\_r6 & 200 & 6& 593 &  79.276 &  75.555 &  68 & 1362.04 &  160.90 \\ 
reg\_n200\_r8 & 200 & 8& 792 &  70.790 &  67.785 &  60 & 1751.27 &  192.74 \\ 
reg\_n200\_r10 & 200&10& 980 &  66.418 &  62.894 &  57 & 2356.94 &  199.12\\
            \hline
        \end{tabular}
    \end{center}
\end{table}

\begin{table}
    \begin{center}
        \caption{Stable set results for graphs from the Erd\H{o}s-R\'enyi model 
        $G(n,p)$}
        \label{randomgraphs}
        \setlength\tabcolsep{5pt}
        \begin{tabular}{|rrr|rrc|rr|} \hline
            graph  &$ n$ & $m$ & $\vartheta(G)$ 
            & \shortstack{final\\ ESB} & $\alpha(G)\geqslant$
            & \shortstack{time\\oracle} & \shortstack{time \\ other} \\ \hline
rand\_n100\_p004 & 100 & 212 & 46.067 &  45.021 &  45 &    338.77&   93.34 \\ 
rand\_n100\_p006 & 100 & 303 & 40.361 &  38.439 &  38 &    769.80&  117.33\\ 
rand\_n100\_p008 & 100 & 443 & 34.847 &  32.579 &  32 &   1126.16&  135.89\\ 
rand\_n100\_p010 & 100 & 489 & 34.020 &  32.191 &  32 &   1004.05&  134.39\\ 
rand\_n200\_p002 & 200 & 407 & 95.778 &  95.032 &  95 &    679.90&  155.80\\ 
rand\_n200\_p003 & 200 & 631 & 83.662 &  81.224 &  80 &   1672.94&  193.75\\ 
rand\_n200\_p004 & 200 & 816 & 73.908 &  70.839 &  67 &   2035.58&  191.42\\ 
rand\_n200\_p005 & 200 & 991 & 69.039 &  66.091 &  62 &   2195.78&  215.50\\  
            \hline
        \end{tabular}
    \end{center}
\end{table}

As a second class of problems we consider random near-$r$-regular graphs, which 
we 
generate as follows. We select a perfect matching on $nr$ vertices and then we
identify consecutive groups of $r$ vertices into a single vertex. This yields
a regular multigraph on $n$ vertices. We remove loops and multiple edges 
resulting in a near-regular graph. 
In Tables~\ref{regulargraphs} and~\ref{randomgraphs} we 
provide results for random graphs.
We compare near-regular graphs with random graphs from the 
Erd\H{o}s-R\'enyi model where the density $p$
is chosen so that the number of edges roughly matches those of the regular 
graphs. We compute our ESB and use a heuristic to compute large stable sets.
In the results the gap between $\vartheta(G)$ and $\alpha(G)$ seems to be 
bigger for regular graphs, but we see in both cases that the ESB reduce the gap 
between $\vartheta(G)$ and
the cardinality of the largest stable set found in a nontrivial way.
Concerning running times we observe the same behavior as before.

\begin{table}
    \setlength\tabcolsep{5pt}
    \centering
    \caption{Tighten $\vartheta(G)$ towards $\alpha(G)$}
    \label{tab: ttf ss}  
    \begin{tabular}{|lrr|rrc|rr|}
        \hline 
        graph & $n$ & $m$ & $\vartheta(G)$  & \shortstack{final\\ ESB} & 
        $\alpha(G)\geqslant$
        & \shortstack{time\\oracle} & \shortstack{time \\ other} \\ \hline    
Circulant47\_030&  47 &  282 & 14.302 & 13.019 & 13 &   193.31&   39.95\\ 
PaleyGraph61    &  61 &  915 &  7.810 &  7.027 &  5 &   319.37&   51.79 \\ 
hamming6\_4     &  64 & 1312 &  5.333 &  4.005 &  4 &    71.37&   13.93 \\ 
spin5           & 125 &  375 & 55.902 & 50.004 & 50 &   342.20&   71.68 \\ 
keller4         & 171 & 5100 & 14.012 & 13.505 & 11 &  1923.60&   62.34\\ 
sanr200\_0\_9   & 200 & 2037 & 49.274 & 47.614 & 42 &  5171.69&  253.94\\ 
c\_fat200\_5    & 200 &11427 & 60.345 & 58.001 & 58 &  3601.56&   85.42\\ 
p\_hat300\_1    & 300 & 8773 &160.345 &160.342 &158 &  5106.69&  112.59\\ 
p\_hat300\_2    & 300 & 6092 &160.345 &158.000 &158 & 31228.51&  187.41\\ 
p\_hat300\_3    & 300 & 3072 &160.345 &158.000 &158 & 19563.05&  252.68\\ 
        \hline 
    \end{tabular} 
\end{table}

As a last experiment for the stable set problem in Table~\ref{tab: ttf ss} we 
consider
instances from the literature, 
taken mostly from the DIMACS challenge~\cite{DIMACS1992}. 
On some instances there is hardly any
improvement of the bound, while other instances are solved to optimality. 
It requires future research to get a better understanding for the fluctuation 
in quality on these instances, 
but for almost all instances the bound improves by at least one integer value.

The computation times for these instances range from 200 to 500 seconds for
the smaller instances ($n \leqslant 125$) to several hours for the biggest 
graphs. As in the instances before a faster oracle would improve the running 
times substantially.

Note that in our computations we aim for getting as good bounds as possible.
If one wants to use the bounds in a branch-and-bound setting, a much more 
aggressive 
strategy with increasing $\kI$ faster and stopping as soon as we do not expect 
to 
reach the next integer value is favorable.

\begin{table}
    \setlength\tabcolsep{5pt}
    \centering
    \caption{Tighten $t^{*}(G)$ towards $\chi(G)$}
    \begin{tabular}{|lrr|rrc|r|}
        \hline 
        graph & $n$ & $m$ & $t^{*}(G)$  & 
        \shortstack{final\\ ESB} & 
        $\chi(G)\leqslant$
        & \shortstack{time\\overall}
        \\ \hline
 myciel3        &  11 &    20 &  2.400 &  3.276 &   4 &     236.95 \\ 
 myciel4        &  23 &    71 &  2.529 &  3.505 &   5 &    1422.71 \\ 
 myciel5        &  47 &   236 &  2.639 &  3.510 &   6 &    4240.61 \\ 
 mug88\_1       &  88 &   146 &  3.000 &  3.022 &   4 &    4709.40 \\ 
 1\_FullIns\_4  &  93 &   593 &  3.124 &  3.939 &   5 &    7219.83 \\ 
 myciel6        &  95 &   755 &  2.734 &  3.534 &   7 &    1540.82 \\ 
 myciel7        & 191 &  2360 &  2.820 &  3.582 &   8 &    2295.24 \\ 
 2\_FullIns\_4  & 212 &  1621 &  4.056 &  4.700 &   6 &   10106.29 \\ 
 flat300\_26\_0 & 300 & 21633 & 16.998 & 17.121 &  26 &    7535.75 \\         
        \hline 
    \end{tabular}
    \label{tab: ttf col}
\end{table}

Results for a selection of 
coloring instances from~\cite{COLInst} are provided  in 
Table~\ref{tab: ttf col}. 
As in the case of 
the stable set problem we use Algorithm~\ref{alg: general procedurce bound 
computation} to obtain ESBs. 
We include at most 100 ESCs in each cycle, 
only include ESCs for subgraphs of order $k \leqslant 8$ 
and perform at most 25 cycles. 
The results are similar in quality to those for stable set from 
Table~\ref{tab: ttf ss},
so for the most instances we are able to obtain bounds, which are one integer 
value better than the original bounds from $t^{*}(G)$. The large running times 
are due to the difficult basic SDP relaxation~\eqref{relaxation col}.

\subsection{The Max-Cut Problem}

Finally we are ready to present computational results for the Max-Cut problem.
It is well known that in the basic SDP relaxation of Max-Cut~\eqref{relaxation 
mc} all ESCs of order 3
can equivalently be represented by the metric polytope~\cite{LaurentPoljak}.
Optimizing over it
gives the exact solution to Max-Cut on graphs not contractible to $K_5$, 
in particular on planar graphs. It is also well known that optimizing over
the metric polytope may lead to rather weak relaxations for general graphs. 
In contrast, the simple SDP relaxation~\eqref{relaxation mc}
provides an upper bound 
at most 14\% above  the optimal value of Max-Cut for graphs with 
nonnegative edge weights, see~\cite{GoemansWilliamson}.

In our computational experiments with Max-Cut we noted that the number of ESCs
necessary to insure that all ESCs for a given value $k$ are satisfied can be
quite large (see Section~\ref{sec:structuralDiff}), even for small values of $n$, such as $n=100$. 
We therefore simplify the ESC relaxation further. If a subgraph $\GI$ violates
the ESC, then instead of asking that $X_I \in CUT_k$, 
we generate a single linear inequality separating $X_I$ from $CUT_k$ and include
it instead of the ESC.  
This weakens the relaxation, but also reduces the computational effort, so that
the total number of ESCs in the model may be quite large, and we can still compute the ESB. 
The computational effort is quite moderate, requiring no more than about 120 
seconds for each of the instances.

\begin{table}
    \setlength\tabcolsep{5pt}
    \centering
        \caption{Max-Cut results for graphs from the Erd\H{o}s-R\'enyi model 
        $G(n,p)$}
        \label{MC-table-1}
        \begin{tabular}{|r|r|rr|r|r|}
            \hline
       $n$ &    $p$       &~~~~~~ 3   & ~~~~~  7  & $z_{mc}$ & \#ESC\\ 
\hline
         &   0.10 & 0.70 & 0.07 & 118 &  2687 \\
     100 &   0.25 & 3.77 & 0.86 & 180 &  3705 \\
         &   0.50 & 5.20 & 2.53 & 246 &  3521 \\
            \hline
         &   0.10  & 7.50 & 4.78 & 184 & 4755 \\
    150  &   0.25  & 7.39 & 5.02 & 310 & 4779 \\
         &   0.50  & 9.71 & 7.51 & 459 & 4605 \\ \hline
        \end{tabular}
\end{table}

We first consider random graphs on $n$ vertices from the Erd\H{o}s-R\'enyi 
model $G(n,p)$.
Each edge is then assigned the weight~$1$  or~$-1$ (each with probability 
$1/2$). 
In Table~\ref{MC-table-1} we report our computational results for $n\in \{100, 
150 
\}$ and 
$p \in \{0.1, 0.25, 0.5 \}$. 
We compare the  ESB with  $k=3$ 
(column labeled 3) to the  ESB with $k = 7$ (column labeled 7). 
The column labeled $3$ provides the deviation 
(in \%) of  
the ESB with $k=3$ from $z_{mc}$. 
Thus if $p$ is the value in the column
labeled $3$, then the ESB is equal to $(1 + p/100)z_{mc}$. 
The column labeled  7 is to be understood in an analogous way
for  $k=7$. 
In all cases we note a substantial gap reduction going from $k=3$ to $k=7$.
The last column contains the number of ESCs at termination. 
It ranges from about 3000 for $n=100$ to about 4500 for $n=150$ and
justifies our strategy to represent each ESC through a single cutting plane.

Next we consider 
graphs from the Beasley 
collection~\cite{BiqMacHomepage} with 
$n=250$. 
Rendl, Rinaldi and 
Wiegele~\cite{RendlRinaldiWiegele} used 10 of these instances 
in a branch-and-bound setting. 
The ``hardest'' instance  250-08 reported in~\cite{RendlRinaldiWiegele} 
resulted in 4553 nodes in the branch-and-bound tree and took 
several days of computation time. 
All the other 9 instances from this
collection resulted in branch-and-bound trees having between 17 and 223 
nodes with computation times in the order
of hours, see Table~6 from~\cite{RendlRinaldiWiegele}.
We recomputed the root bound for all these instances and present our root gap 
in Table~\ref{MC-literature}.
We find it remarkable that our new bounding procedure is strong enough
to prove optimality for all these instances right at the root node 
with the exception of problem 250-08.
For this problem the gap at the root node was 2.19\%. We recomputed the root 
bound 
in our setting and came up with a root gap of only 0.5\%, thus
reducing the gap by 75\%.

\begin{table}
    \setlength\tabcolsep{5pt}
    \centering
        \caption{Max-Cut results for graphs from the OR library}
        \label{MC-literature}
        \begin{tabular}{|r|r|rr|r|}
            \hline
            graph   & opt.\ cut   & \# branch-and-bound 
            nodes~\cite{RendlRinaldiWiegele} & root 
            gap~\cite{RendlRinaldiWiegele} & our root gap\\ 
\hline
250-01 & 45607   &   37  & 0.44& * \\
250-02 & 44810   &   19  & 0.56& * \\
250-03 & 49037   &   19  & 0.14& * \\
250-04 & 41274   &   17  & 0.39& * \\
250-05 & 47961   &   21  & 0.35& * \\ 
250-06 & 41014   &  223  & 1.03& * \\
250-07 & 46757   &   37  & 0.44& * \\
250-08 & 35726   & 4553  & 2.19& 0.5 \\
250-09 & 48916   &   47  & 0.78& * \\
250-10 & 40442   &   63  & 0.62& * \\ 
    
            \hline
        \end{tabular}
\end{table}

\begin{table}
    \setlength\tabcolsep{5pt}
    \centering
        \caption{Max-Cut results for Chimera graphs with $n=512$}
        \label{MC-chimera}
        \begin{tabular}{|r|rc|r|}
            \hline
            graph       & final ESB   & best found cut   & \#ESC\\ 
\hline
chimera-1 & 434.38   & 433  & 22275 \\
chimera-2 & 452.69   & 451  & 24707 \\
chimera-3 & 447.35   & 447  & 22390 \\
chimera-4 & 439.90   & 439  & 20748 \\
chimera-5 & 441.66   & 440  & 22838 \\      
            \hline
        \end{tabular}
\end{table}

As a final experiment we consider Max-Cut instances on 
{\em Chimera graphs}. This class of graphs has found increased interest 
in connection with quantum annealing, see~\cite{Chimera} for further details. 
In Table~\ref{MC-chimera} we
provide computational results with such graphs on $n=512$ vertices. We compute 
our ESB and also use a heuristic to find a good cut.
It turns out that our bounding approach works nicely on these graphs,
leading to provably optimal solutions in 2 out of 5 instances and the smallest 
possible positive gap (of 1) in the remaining cases. 
The computation times for each of these (big) instances range from
700 to 900 seconds, which we consider remarkable when dealing with more than 
20000 ESCs.

We conclude that for Max-Cut 
our ESB constitute a substantial
improvement compared to the previously used strongest bounds based
on SDP with triangle inequalities~\cite{RendlRinaldiWiegele}. These correspond to the column 
labeled $3$ in Table~\ref{MC-table-1}.

\section{Conclusions and Future Work}
\label{sec:conclusions}

Summarizing, we offer the following conclusions from the 
computational results.
Our computational approach based on the partial Lagrangian dual is very efficient in handling also a large number of ESCs. The dual function 
evaluation separates the SDP part from the ESCs and 
therefore opens the way for large-scale computations. The minimization 
of the dual function is carried out as a convex quadratic optimization problem without any SDP constraints, and therefore is also suitable for a large number of 
ESCs. 

Our computational results for 
stable set and coloring confirm the theoretical hardness results
for these problems. Including ESCs of rather small size
($k \leqslant 8$) yields a noticeable improvement of the bounds.

The limiting factor for stable set instances is the solution time 
of the oracle. 
Hence it is desirable to have a fast 
solver for these kind of instances.

On the practical side we consider the 
cutting plane weakening of the ESCs 
for Max-Cut a promising new way to tighten bounds for this problem. 

It will be a future project to explore these bounds in a branch-and-bound 
setting in order to solve Max-Cut, stable set and coloring 
instances to optimality.

%
%

 \bibliographystyle{splncs04}
 \bibliography{papers}

\begin{thebibliography}{10}
\providecommand{\url}[1]{\texttt{#1}}
\providecommand{\urlprefix}{URL }
\providecommand{\doi}[1]{https://doi.org/#1}

\bibitem{AARW}
Adams, E., Anjos, M.F., Rendl, F., Wiegele, A.: A {H}ierarchy of subgraph
  projection-based semidefinite relaxations for some {NP}-hard graph
  optimization problems. INFOR Inf. Syst. Oper. Res.  \textbf{53}(1),  40--47
  (2015)

\bibitem{IntroductionSOCP}
Alizadeh, F., Goldfarb, D.: Second-order cone programming. Math. Program.
  \textbf{95}(1, Ser. B),  3--51 (2003), iSMP 2000, Part 3 (Atlanta, GA)

\bibitem{BiqMacHomepage}
{Biq Mac Library}: \url{http://biqmac.aau.at/}, {Last accessed 15 May 2019}

\bibitem{NumericalOptimizationBonnans}
Bonnans, J.F., Gilbert, J.C., Lemar{\'e}chal, C., Sagastiz\'{a}bal, C.A.:
  Numerical Optimization: Theoretical and Practical Aspects. Springer-Verlag,
  Secaucus, NJ, USA (2006)

\bibitem{BorosCramaHammer}
Boros, E., Crama, Y., Hammer, P.L.: Upper-bounds for quadratic {$0$}-{$1$}
  maximization. Oper. Res. Lett.  \textbf{9}(2),  73--79 (1990)

\bibitem{Chimera}
Dash, S., Puget, J.F.: On quadratic unconstrained binary optimization problems
  defined on {C}himera graphs. Optima  \textbf{98} (2015)

\bibitem{DADAL}
{De Santis}, M., {Rendl}, F., {Wiegele}, A.: {Using a Factored Dual in
  Augmented Lagrangian Methods for Semidefinite Programming}. ArXiv e-prints
  (Oct 2017)

\bibitem{DelormePoljak}
Delorme, C., Poljak, S.: Laplacian eigenvalues and the maximum cut problem.
  Math. Programming  \textbf{62}(3, Ser. A),  557--574 (1993)

\bibitem{DIMACS1992}
{DIMACS Implementation Challenges}: \url{http://dimacs.rutgers.edu/Challenges/}
  (1992), {Last accessed 15 May 2019}

\bibitem{FischerGruberRendlSotirovBundleMaxCutEquipartition}
Fischer, I., Gruber, G., Rendl, F., Sotirov, R.: Computational experience with
  a bundle approach for semidefinite cutting plane relaxations of {M}ax-{C}ut
  and equipartition. Math. Program.  \textbf{105}(2-3, Ser. B),  451--469
  (2006)

\bibitem{EasyBundle}
Frangioni, A., Gorgone, E.: Bundle methods for sum-functions with ``easy''
  components: applications to multicommodity network design. Mathematical
  Programming  \textbf{145}(1),  133--161 (2014)

\bibitem{elli-diss}
Gaar, E.: Efficient Implementation of SDP Relaxations for the Stable Set
  Problem. Ph.D. thesis, Alpen-Adria-Universit\"at Klagenfurt (2018)

\bibitem{GaarRendl}
Gaar, E., Rendl, F.: A bundle approach for {SDP}s with exact subgraph
  constraints. In: Lodi, A., Nagarajan, V. (eds.) Integer Programming and
  Combinatorial Optimization. pp. 205--218. Springer International Publishing
  (2019)

\bibitem{GoemansWilliamson}
Goemans, M.X., Williamson, D.P.: Improved approximation algorithms for maximum
  cut and satisfiability problems using semidefinite programming. J. Assoc.
  Comput. Mach.  \textbf{42}(6),  1115--1145 (1995)

\bibitem{OurUsedFormOfLovasTheta}
Gr\"otschel, M., Lov\'asz, L., Schrijver, A.: Geometric algorithms and
  combinatorial optimization, Algorithms and Combinatorics: Study and Research
  Texts, vol.~2. Springer-Verlag, Berlin (1988)

\bibitem{GuruswamiKhanna}
Guruswami, V., Khanna, S.: On the hardness of 4-coloring a 3-colorable graph.
  SIAM J. Discrete Math.  \textbf{18}(1),  30--40 (2004)

\bibitem{stableSetNotApproximable}
H{\aa}stad, J.: Clique is hard to approximate within {$n^{1-\epsilon}$}. Acta
  Math.  \textbf{182}(1),  105--142 (1999)

\bibitem{HelmbergRendlBundleSDP}
Helmberg, C., Rendl, F.: A spectral bundle method for semidefinite programming.
  SIAM J. Optim.  \textbf{10}(3),  673--696 (2000)

\bibitem{HUL2}
Hiriart-Urruty, J.B., Lemar{\'e}chal, C.: Convex Analysis and Minimization
  Algorithms {II}: Advanced Theory and Bundle Methods. Grundlehren der
  mathematischen Wissenschaften, Springer-Verlag (1993)

\bibitem{Khot}
Khot, S.: Improved inapproximability results for {M}ax{C}lique, chromatic
  number and approximate graph coloring. In: 42nd {IEEE} {S}ymp. on {F}ound. of
  {C}omp. {S}c. ({L}as {V}egas, {NV}, 2001), pp. 600--609. IEEE Comp. Soc., Los
  Alamitos, CA (2001)

\bibitem{KiwielProximityControl}
Kiwiel, K.C.: Proximity control in bundle methods for convex nondifferentiable
  minimization. Mathematical Programming  \textbf{46}(1),  105--122 (Jan 1990)

\bibitem{LasserreHierarchy}
Lasserre, J.B.: An explicit exact {SDP} relaxation for nonlinear 0-1 programs.
  In: Integer programming and combinatorial optimization ({U}trecht, 2001),
  Lecture Notes in Comput. Sci., vol.~2081, pp. 293--303. Springer, Berlin
  (2001)

\bibitem{LaurentPoljak}
Laurent, M., Poljak, S.: The metric polytope. In: Balas, E., Cornuejols, G.,
  Kannan, R. (eds.) Integer Programming and Combinatorial Optimization. pp.
  247--286 (1992)

\bibitem{LovaszStart}
Lov{\'a}sz, L.: On the shannon capacity of a graph. IEEE Transactions on
  Information Theory  \textbf{25}(1), ~1--7 (1979)

\bibitem{LovaszSchrijverHierarchy}
Lov\'asz, L., Schrijver, A.: Cones of matrices and set-functions and
  {$0$}-{$1$} optimization. SIAM J. Optim.  \textbf{1}(2),  166--190 (1991)

\bibitem{mosek}
{MOSEK ApS}: The MOSEK optimization toolbox for MATLAB manual. Version 8.0.
  (2017), \url{http://docs.mosek.com/8.0/toolbox/index.html}

\bibitem{COLInst}
Nguyen, T.H., Bui, T.: Graph coloring benchmark instances.
  \url{https://turing.cs.hbg.psu.edu/txn131/graphcoloring.html}, {Last accessed
  15 May 2019}

\bibitem{RendlRinaldiWiegele}
Rendl, F., Rinaldi, G., Wiegele, A.: Solving max-cut to optimality by
  intersecting semidefinite and polyhedral relaxations. Math. Program.
  \textbf{121}(2, Ser. A),  307--335 (2010)

\bibitem{RendlSotirovBundleSDPQAP}
Rendl, F., Sotirov, R.: Bounds for the quadratic assignment problem using the
  bundle method. Math. Program.  \textbf{109}(2-3, Ser. B),  505--524 (2007)

\bibitem{RobinsonBundlePrimalConvergence}
Robinson, S.M.: Bundle-based decomposition: conditions for convergence. Annales
  de l'I.H.P. Analyse non lin\'eaire  \textbf{S6},  435--447 (1989)

\bibitem{RockafellarConvexAnalysis}
Rockafellar, R.T.: Convex Analysis. Princeton Mathematical Series, No. 28,
  Princeton University Press, Princeton, N.J. (1970)

\bibitem{SheraliAdamsHierarchy}
Sherali, H.D., Adams, W.P.: A hierarchy of relaxations between the continuous
  and convex hull representations for zero-one programming problems. SIAM J.
  Discrete Math.  \textbf{3}(3),  411--430 (1990)

\bibitem{SDPT31}
Toh, K.C., Todd, M.J., T\"{u}t\"{u}nc\"{u}, R.H.: S{DPT}3---a {MATLAB} software
  package for semidefinite programming, version 1.3. Optim. Methods Softw.
  \textbf{11/12}(1-4),  545--581 (1999)

\bibitem{SDPT32}
T\"{u}t\"{u}nc\"{u}, R.H., Toh, K.C., Todd, M.J.: Solving
  semidefinite-quadratic-linear programs using {SDPT}3. Math. Program.
  \textbf{95}(2, Ser. B),  189--217 (2003)

\bibitem{SDPNALPlus}
Yang, L., Sun, D., Toh, K.C.: {${\rm SDPNAL}+$}: a majorized semismooth
  {N}ewton-{CG} augmented {L}agrangian method for semidefinite programming with
  nonnegative constraints. Math. Program. Comput.  \textbf{7}(3),  331--366
  (2015)

\end{thebibliography}

\end{document}